
\input amstex

\font \MySmallFont= cmr8

\newcount\secnum \secnum=0       
\newcount\subsecnum              

\def\section#1{\advance\secnum by 1 \subsecnum=0%
            \head{\the\secnum. #1}\endhead }

\def\subsection#1{\advance\subsecnum by 1%
       \subhead\nofrills
       {\the\subsecnum. #1.\ }\endsubhead}
\newcount\firstpageno
\newcount\tocpageno

\def\d{{\bold d}}
\def\e{{\bold e}}

\def\r{{\bold r}}
\def\s{{\bold s}}

\def\u{{\bold u}}
\def\v{{\bold v}}
\def\w{{\bold w}}

\def\B{{\bold B}}

\def\L{{\bold L}}
\def\N{{\bold N}}

\def\R{{\bold R}}
\def\S{{\bold S}}
\def\Reals{{\Bbb R}}

\def\num{\bold n}

\def \ma(#1;#2;#3;#4;#5;#6;#7;#8;#9)%
{\left(\matrix#1&#2&#3\\ #4&#5&#6\\ #7&#8&#9\endmatrix\right)}

\def \ve(#1;#2;#3)%
{\langle\matrix#1, &#2, &#3\endmatrix\rangle}

\def \ip(#1,#2)%
   {\langle #1,#2\rangle}

\def\today{\ifcase\month\or}

\def \Im{\hbox{Im}}

\def \Jac{\hbox{Jac}}

\def\Aff{{\bold {Aff}}}

\def \mutilde {\tilde \mu}

\def \Prob(#1){\hbox{Prob}(#1)} 

\def \Rnd{\hbox{Rnd}}

\def\Diff{\hbox{Diff}}

\documentstyle{amsppt}

\refstyle{A}

\nologo

\def\proof {\noindent{\it Proof.}\enspace}
\def\qedmark{\hbox{\vrule height 4pt width 3pt}}
\def\qedskip{\vrule height 4pt width 0pt depth 1pc}
\def\qed{\nobreak\quad\nobreak{\qedmark\qedskip}}
\def\eps{\epsilon}

\vsize = 9in
\hsize = 6in
\hoffset = 20pt

\topmatter

\title{Point Clouds: Distributing Points \break Uniformly on A Surface}\endtitle

\author Richard Palais, Bob Palais, and Hermann Karcher\endauthor



\endtopmatter

\document

{\narrower\narrower\narrower \overfullrule=0pt
\noindent
{\bf ABSTRACT: }
   
         The concept of ``Point Cloud'' has played an increasingly important  r\^ole in many areas of 
Engineering, Science, and Mathematics.  Examples are: LIDAR,  3D-Printing, Data 
Analysis, Computer Graphics, Machine Learning, Mathematical Visualization,  
Numerical Analysis, and Monte Carlo Methods.  Entering ``point cloud" into Google 
returns nearly 3.5 million results!   A {\it point cloud} for a finite volume manifold $M$ 
is a finite subset or a sequence in $M$, with the essential feature that it is a 
``representative" sample of $M$. The definition of a point cloud varies with its use, 
particularly what constitutes being representative.  Point clouds arise in many different ways:
in LIDAR they are just 3D data captured by a scanning device, while in Monte 
Carlo applications they are constructed using highly complex algorithms developed 
over many years. In this article we outline a rigorous mathematical theory of point 
clouds, based on the classic Cauchy-Crofton formula of Integral Geometry and its 
generalizations. We begin with point clouds on surfaces in $\R^3$, which simplifies 
the exposition and makes our constructions easily visualizable. We proceed to hypersurfaces 
 and then submanifolds of arbitrary codimension in $\R^n$, and finally, using an elegant 
 result of J\"urgen Moser, to arbitrary smooth manifolds with a volume element.   \par

\section{Introduction}
   While the mathematical ideas to be discussed in this article have a theoretical flavor, they 
in fact evolved from our need to overcome vexing practical problems that arose during our   
development of mathematical visualization software.
\footnote{3D-XplorMath; see http://3D-XplorMath.org. Also, see [P].}
Those problems centered on difficulties we encountered while trying to devise appropriate 
rendering methods for displaying surfaces that are implicitly defined. Important examples of such 
surfaces are often of large genus, being defined as level sets of polynomials of high 
degree, and as a result cannot be explicitly parameterized. In addition they may 
be immersed rather than embedded and have singularities and self-intersections that both 
add to their interest but also make them difficult to visualize. A major problem contributing 
to the difficulty of working algorithmically with implicitly defined surfaces is that it is much harder 
to ``choose'' a point on an implicit surfaces than on one given parametrically, a consideration 
that we will discuss in more detail later. These inter-related factors make the customary approach 
to rendering implicit surfaces, based  on ``raytracing'',  much slower and less flexible than the 
various methods available for rendering parametric surfaces, and also make it more difficult to display 
their important structural features. Not surprisingly then, we have kept searching for improved 
methods for visualizing implicit  surfaces.\par

   In raytracing a surface, one carries out a complex process for each pixel on the screen, even for those 
that turn out not to show a point of the surface. Moreover, even as graphics hardware has improved in 
speed, screen resolutions have increased in step, so that raytracing remains a time-consuming process. 
About twelve years ago it occurred to us to try a ``Pointilist'' approach---that is to represent our implicit surfaces 
by  ``point clouds'', i.e.,  relatively dense sets of points distributed evenly over the surface. We guessed that 
if we could create such clouds of points and render them in anaglyph stereo (the kind requiring red/green 
glasses), then not only could we get away with processing far fewer pixels---for a major speedup---but as a 
side-effect we would automatically achieve the transparency required to show all parts of a surface simultaneously, 
and in this way reveal otherwise hidden structural details that we often wished to display. Early experiments were highly 
encouraging, confirming these hopes and expectations.  In particular, when rotating a surface, a point cloud, 
unlike  a raytrace, only needs to be computed once, so even when using slower hardware it gave 
an important speed-up.  But a major problem remained before these experiments could mature into a fully 
successful rendering technique;  namely we had to develop effective and efficient algorithms for generating 
point clouds that had the required uniform density. 

  In the remainder of this article we will explain the mathematical ideas behind our solution to that problem
and its generalizations and some interesting explorations we were led to while working on it.
Using few (and metaphorical) words, what eventually worked 
was throwing random darts and seeing where they hit the surface, or in less picturesque language, 
seeing where randomly selected lines intersected the surface. This approach was both motivated and 
justified by the famous Cauchy-Crofton Formula from Integral Geometry, which says that the 
area of a piece of surface is proportional to the average number of times lines intersect it. 

  But a {\it lot\/} is hidden in that word ``random''!  Our use of the term has little relation to the randomness in 
physical processes like radioactivity, or the way repeated measurements with a physical device are distributed 
about their mean.  Rather, as our experiments continued, we quickly discovered that to make our point clouds 
look  ``right'',  it was essential for the points of the cloud to be ``equidistributed with respect to area'',  meaning 
roughly that parts of the surface with approximately equal areas should contain approximately equal numbers 
of cloud points. Finding out how to make rigorous mathematical sense of that, and how to implement it algorithmically 
was a key step in developing our point cloud rendering techniques. It  led us on a learning adventure into an old 
and fascinating area of mathematics;  namely how to use a Random Number Generators (RNG) to distribute 
points uniformly in all manner of interesting spaces, and we hope that the results presented below will count 
as a useful contribution to that circle of ideas. 

\section {Review of Required Background}

\subsection {Some Fixed Notation, Definitions, and Terminology}

   We will denote by $\B^n$ and $\S^{n-1}$ respectively the open unit ball and the unit sphere in $\Reals^n$
 (centered at the origin), and for $v\in \S^{n-1}$ we denote by $v^\perp$  the $(n-1)$-dimensional 
 subspace of $\Reals^n$ orthogonal to $v$. If $D$ is a bounded open subset of $v^\perp$ and $f$ is a
 smooth real-valued  function defined on $D$, then we call the subset $\{ (x,f(x)v)) \mid x \in D\}$
 a {\it smooth graph hypersurface\/} in $\Reals^n$. If $D$ has a (piecewise) smooth boundary, and 
 if $f$ extends to be smooth on the boundary of $D$,  then we call $\{ (x,f(x)v)) \mid x \in \bar D\}$ 
a {\it smooth graph hypersurface with (piecewise) smooth boundary\/}.  A subset $\Sigma$ of $\Reals^n$ 
will be called a  {\it hypersurface in $\Reals^n$\/} if it can be represented as a finite union of smooth 
graph hypersurfaces in $\Reals^n$, some possibly with (piecewise) smooth boundary. 
If $F:\Reals^n \to \Reals$ is a smooth function, then $p \in \Reals^n$ is a {\it critical point\/} of $F$ if $DF_p =0$, 
and $r \in \Reals$ is a {\it regular value\/} of $F$ if the the level set $\Sigma := \{x \in \Reals^n \mid F(x) =r\}$ contains 
no critical points of $F$; it then follows  from the Implicit Function Theorem that any bounded part of such a
$\Sigma$ is a regular hypersurface, and such a hypersurface is what we will mean by an {\it implicit surface\/}.

\subsection{Sequences}

  Sequences arise frequently in what follows and it will be  convenient to establish notation to help deal 
efficiently with our uses of them.  An infinite sequence of points, $x_1, x_2, \ldots$, will be denoted by a symbol 
such as $\{x_i\}$, in which case,  for a positive integer $N$,  $\{x_{N+i}\}$ denotes the infinite subsequence 
$x_{N+1}, x_{N+2}\ldots$, while $\{x_i\}_N$ will denote the complementary  finite sequence $x_1, x_2, \ldots x_N$.
 If $E$ is a set,  $\#(\{x_i\}_N, E)$ will denote the number of indices $i$ with $1\le i \le N$ such that $x_i \in E$. 
Then a natural interpretation of the quotient, $\#(\{x_i\}_N, E)\over N$ (and one appropriate for our purposes) 
is the probability that if $1\le i \le N$ then $x_i \in E$.  This leads us to the following: 

\definition{Definition}  
Let $X$ be a set and  $\{x_i\}$ an infinite sequence in $X$.  For each subset $E \subseteq X$,
if  $\lim_{N \to \infty}{ \#(\{x_i\}_N, E)\over N}$ exists  and equals $L$, then we say that
 $\Prob(\{x_i\} \in E) $ is defined and is equal to $L$,  
 otherwise we say  that $\Prob(\{x_i\} \in E) $ is undefined.
 \enddefinition
 \par

\subsection {Densities and Volume Forms} 

  We will assume familiarity with the theory of volume forms on manifolds and the associated measure 
theory. However these  play such an important role in what follows that we  provide next a quick review 
to establish consistent notation and terminology. For details see Chapter 3 of [Sp-1] and Chapter 7 of [Sp-2].
\par

   In what follows, $M$ denotes a smooth (i.e., $C^\infty$) $n$-dimensional manifold, possibly with a piecewise 
smooth boundary, $\partial M$.  (Much of the time, $M$ will be assumed compact.) $\Diff(M)$ will denote the 
group of diffeomorphisms of $M$ with the $C^\infty$ topology.  A {\it density \/} for $M$,  
$d\mu$, is the absolute value $|\omega|$ of a smooth $n$-form $\omega$, and if $\omega$ does not vanish 
anywhere we call $d\mu$ a {\it volume element\/}. Since the space $\Lambda^n(T^*M_p)$ of $n$-forms at a 
point $p$ of $M$ is one-dimensional, if $d\mu$ is a volume element then any density $d\gamma$ for $M$ can be 
written uniquely in the form $d\gamma = \rho \, d\mu$ where $\rho: M \to \Reals^+$ is a nonnegative smooth function. 
\par 

 In particular, suppose that   $(x_1,\ldots,x_n)$ is a local coordinate system for an open set $O \subseteq M$, 
 given by a diffeomorphism $x: O \to x(O) \subseteq \Reals^n$. 
 Then in $O$ we have the coordinate representations: 
 $\omega = f(x_1,\ldots,x_n)\, |dx_1\wedge \ldots \wedge dx_n|$ 
 and  $d \mu =  |f(x_1,\ldots,x_n)|\, |dx_1\wedge \ldots \wedge dx_n|$ 
 for some smooth function $f: x(O) \to \Reals$.  
If $(\tilde x_1,\ldots,\tilde x_n)$ is a second coordinate system in $O$, 
then the Change of Variable Formula of multi-variable calculus says that the corresponding 
coordinate representations of $\omega$ and  $d \mu$  with respect to these new 
coordinates is given by replacing the function $f(x_1,\ldots,x_n)$ by the function 
$\tilde f(\tilde x_1,\ldots,\tilde x_n) = \Jac(\tilde x,x)\, f(x_1(\tilde x),\ldots,x_n(\tilde x) )$,
where $\Jac(\tilde x,x) := \det \left({\partial \tilde x_i \over \partial x_j}\right)$ is the Jacobian 
determinant of the transformation of coordinates $\tilde x_i = \tilde x_i(x)$ from $x(O)$ to $\tilde x(O)$
and $x_i = x_i(\tilde x)$ is the inverse transformation from $\tilde x(O)$ to $x(O)$.  
We will call the coordinate system  $(x_1,\ldots,x_n)$ {\it canonical coordinates\/} with respect to 
$d\mu$ if $f$ is identically one in $O$, i.e.,   if $\omega = dx_1\wedge\ldots \wedge dx_n$ in $O$ 
(so that $d\mu = |dx_1\wedge\ldots\wedge dx_n|$ in $O$).  For such coordinates, the measure $\mu(K)$ 
of a compact subset $K \subseteq O$ is just the Lebesgue measure of the image $x(K) \subseteq \Reals^n$
under the coordinate map $x: O \to \Reals^n$. 
It is easy to see that at every point $p$ of $M$ we can choose canonical coordinates.  For if $(y_1,\ldots,y_n)$ is 
any coordinate system centered at $p$ and $\omega = f(y_1,\ldots , y_n) \, dy_1\wedge \ldots\wedge dy_n$ then we can
define new coordinates $(x_1,\ldots,x_n)$ by the transformation law 
 $x_1(y_1,\ldots , y_n) := \int_0^{y_1} f(t,y_2,\ldots,y_n) \,dt$ and  $x_i := y_i$ for $i >1$. 
Then $dx_1 = f(y) dy_1 + \sum_{j=2}^n a_j(y) dy_j$ so $dx_1 \wedge(dy_2 \wedge \ldots \wedge dy_n) = \omega$, 
proving that $(x_1,\ldots,x_n)$ are canonical coordinates. 

   The density $d\mu$ determines a Radon measure $\mu$ on $M$. If $g$ is a continuous real-valued 
function on $M$ with support a compact subset of the coordinate domain $O$ above, then 
$\int_M g(p) \, d\mu(p) = \int_O g(x_1,\ldots,x_n) |f(x_1,\ldots,x_n)| dx_1\,\ldots \, dx_n$. 
The integral for a more general $g$ can then be computed using a partition of unity subordinate 
to a covering by coordinate neighborhoods. \par

     If $F:M \to N$ is a diffeomorphism between two compact $n$-dimensional  manifolds and $d\nu$ is a density for
$N$, then its pull-back, $d\mu := DF^*(d\nu)$, is a density for $M$, and if $d\nu$ is a volume element then so is 
$d\mu$.  The corresponding measures, $\nu$ and $\mu$ on $N$ and $M$, are clearly related by  $\nu = \mu \circ F^{-1}$, 
so in particular $\nu(N) = \mu \circ F^{-1}(N) = \mu(M)$, i.e., they have the same total volume, and it follows in particular 
that the natural action of $\Diff(M)$ on the space of volume forms on $M$, given by pull-back, preserves total volumes.


\subsection {Measured Manifolds} 

By a {\it measured manifold\/} we mean a smooth manifold $M$ together with a choice of volume element, 
$d\mu_M$,  and by canonical coordinates for $M$ we will mean canonical coordinates with respect to $d\mu_M$.
Any smooth density $d\gamma$ on $M$ can be written uniquely as  $d\gamma = \alpha \, d\mu_M$ where $\alpha$ 
is a non-negative smooth function on $M$ that  we will refer to as the Radon-Nikodym  derivative of $d\gamma$.  \par

 Various kinds of structured manifolds have natural choices of volume elements. 
For example if $M$ is a symplectic manifold of dimension $n=2k$ 
and $\lambda$ is its symplectic form, then $d\mu_M := |\lambda^k|$ defines $M$ as a 
measured manifold, and similarly, for a Riemannian manifold $M$ the standard choice 
is $d\mu_M := |\omega|$, where for any point $p$ of $M$, $\omega_p = \theta_1 \wedge\ldots \wedge \theta_n$,
with $\theta_1,\ldots,\theta_n$ any orthonormal basis for $T^*M_p$.
Equivalently, in local  coordinates $(x_1,\ldots,x_n)$ , 
$d\mu_M = \sqrt{g}\, dx_1\, dx_2\,\ldots dx_n$, where $g := \det(g_{ij})$ is the determinant 
of the metric tensor. In particular we will consider any smooth submanifold of $\Reals^n$ to 
be a Riemannian,  and hence measured,  manifold, with its induced Riemannian metric. 
More particularly, a linear subspace, $V$,  of $\Reals^n$ is a measured manifold with the associated 
measure being Lebesgue measure, and in this case we will as usual write integrals as $\int_V g(x) \,dx$.

\definition{Definition} If $M$ is a measured manifold, then an embedding $F: M \to \Reals^n$ of $M$  into 
a Euclidean space will be called {\it measure preserving\/} if it pulls back the Riemannian volume form on $F(M)$
to the volume form $d\mu_M$ on $M$.
\enddefinition
      \par 
      
      Let $F:M \to N$ be a smooth map between two measured manifolds of the same dimension.
The pull-back $DF^*(d\mu_N)$ is a smooth density on $M$ and so, as above, it can be 
written in the form $\alpha \, d\mu_M$ for a unique $\alpha: M \to \Reals^+$. 
For the case $M =N = \Reals^n$,  the Change of Variable Theorem says that 
$\alpha = \left |\det \big({\partial F_i\over \partial x_j}\big)\right|$, the absolute 
value of the Jacobian determinant. If we choose canonical coordinates for 
$M$ at $p$ and for $N$ at $F(p)$, then using these coordinates we can consider 
$F$ as a map of $\Reals^n$ to $\Reals^n$, and again $\alpha$ will be the 
absolute value of the Jacobian determinant. For this reason we refer to $\alpha$ as the 
{\it generalized Jacobian} of $F$ and write it as $|\Jac(F)|$. 
Clearly $|\Jac(F)|$ is everywhere non-negative and vanishes precisely on the set $\Xi$  
of critical points of $F$, i.e., the points where $DF_p$ is not invertible, 
or equivalently the points where $F$ is not a local diffeomorphism. The subset $F(\Xi)$ 
of $N$ is called the set of critical values of $F$ and it is an important fact, called 
The Morse-Sard Theorem ([H] page 69 or [SP-1] p. 72) that it has measure zero.  
With this formalism, the change of variable theorem takes the following form. 

\proclaim{Proposition. (Change of Variable Formula for Measured Manifolds)}
Let $F:M \to N$ be a smooth injective map between two measured manifolds of the same dimension,  
and let $g$ be an integrable  function on $N$. Then  $g\circ F$ is integrable on $M$ and 
$$\int_{N} g(y) \, d\mu_{N}(y)  =  \int_{M} g(F(x))\, |\Jac(F)|(x) \, d\mu_{M}(x).$$
In particular $\mu(N) = \int_{M} |\Jac(F)| \, d\mu_{M}$.
\endproclaim

\subsection {A Theorem of J\"urgen Moser}

   As was noted above, for a compact, smooth manifold $M$, the natural action (given by pull-back) 
of $\Diff(M)$ on the space of volume forms on $M$ preserves total volumes.  In [M], Moser gave a 
straightforward and elegant proof of the following ``converse'' result. 

\proclaim{Moser's Theorem} If $M$ is a compact, smooth manifold then $\Diff(M)$ acts transitively on the 
space of volume forms on $M$ having any fixed total volume. In particular, if $M$ is a compact measured 
manifold and $d\nu$ is any volume element for $M$ with the same total volume as $d\mu_M$, then 
there is an $\phi \in \Diff(M)$ such that $d\mu_M = D\phi^*(d\nu)$.
\endproclaim

\proclaim{Corollary (Moser-Whitney Theorem)}
A compact, smooth, measured manifold $M$ of dimension $n$ admits a volume preserving embedding 
in $R^{2n}$.
\endproclaim

\proof By a theorem of Hassler Whitney [Wh], there exists a smooth embedding $F$ of $M$ onto 
a submanifold $M'$ of $R^{2n}$.  Scaling $R^{2n}$ by a suitable positive constant $c$ will 
map $M'$ diffeomorphically onto a submanifold $M''$ of $R^{2n}$ with the same volume as $M$, 
hence, replacing $F$ by $cF$, we can assume that the total volume of $M$ and of $M'$ are equal.
Let $\nu$ denote the (Riemannian) volume element of $M'$. Then, since $F$ is a diffeomorphism, 
$d\mu_M$ and $DF^*(\nu)$ are volume elements for $M$ with the same total volume and, by Moser's 
Theorem, there exists a diffeomorphism $\phi \in \Diff(M)$ such that $d\mu_M = D\phi^*(DF^*(\nu))$, 
so that $F\circ\phi: M \to \Reals^{2n}$ gives the desired embedding.   \qed
\par

\subsection {Maps Between Measured Manifolds and The Co-area Formula}

  We will need the following important generalization of the Change of Variable Formula that handles 
the case that $F$ is not assumed to be injective.  We write $\#(X)$ for the number of points in a set $X$.

\proclaim{Theorem (Co-area Formula for measured manifolds of the same dimension)} \par
If $F$ is a smooth map from a compact measured manifold $M$ to another smooth manifold $N$ 
of the same dimension then   $$ \int_N \#(F^{-1}(p)) \, d\mu_N(p) = \int_M |\Jac(F)| \, d\mu_M.$$
\endproclaim

\proof  For any open subset $O$ of $M$  we will denote by $F_O: O \to N$ the restriction of $F$ to $O$ 
and we will prove that the formula $ \int_N \#({F_O}^{-1}(p)) \, d\mu_N(p) = \int_O |\Jac(F)| \, d\mu_M$ 
holds for all $O$. Of course the Co-Area Formula is just the special case where we take $O = M$, 
but using the more general statement will give us extra flexibility that we will make use of in the proof. \par

  We first consider the case that $F_O$ is injective. Then $\#({F_O}^{-1}(p))$ is identically one on $F(O)$ 
and zero on its complement, i.e., it is the characteristic function $\chi(p)$ of $F(O)$, and the claimed 
formula reduces to $\mu_N(F(O)) = \int_O |\Jac(F)| \, d\mu_M$. So this first case follows from the 
preceding proposition by taking $F$ to be $F_O$ and $g$ to be  $\chi$. 
\par

The set $\Xi$  of critical points of $F$ is clearly closed in $M$, and since $\Jac(F)$ vanishes on 
$\Xi$, and $F(\Xi)$ has measure zero, we can replace the open set $O$ by $O \setminus \Xi$ 
without changing either side of the claimed formula.  So, without loss of generality, we can 
assume that $F_O$ has no critical points and therefore, by the implicit function theorem, that it is 
a local diffeomorphism. Thus, for each $p\in O$ we can choose an open neighborhood $U$ of $p$ 
for which $F_U$ is a diffeomorphism onto $F(U)$, and we can choose $U$ so small that its 
closure $\bar U$ is compact and its boundary, $\partial U := \bar U \setminus U$ has measure zero.
Since these $U$ form an open covering of $O$,  by Lindel\"of's Theorem we can 
find a sequence $\{U_i\}$ of these $U$ that already cover $O$.   By choice of the $U_i$, 
each $F_{U_i}$ is a diffeomorphism onto its image, so by the first case we have 
$ \int_N \#({F_{U_i}}^{-1}(p)) \, d\mu_N(p) = \int_{U_i} |\Jac(F)| \, d\mu_M$, and moreover 
the same holds if we replace $U_i$ by any open set  $O_i \subseteq U_i$.

\par

  We next observe that for each $p$ in $N$, the function $V \mapsto   \#({F_V^{-1}(p))}$ is a countably 
additive function on the subsets of $N$;  in fact it is the measure given by the sum of delta measures one at 
each point of $F^{-1}(p)$.  To use this observation we ``disjointify'' the sequence $\{U_i\}$, that is 
we replace it by the sequence $O_i$ of disjoint open sets given by $O_1 := U_1$ and 
$O_{i+1} := U_i \setminus \{\overline {U_1} \cup \overline {U_2} \cup \ldots \cup \overline {U_i}\}$.
While these $O_i$ no longer cover $O$, their union $O^\prime := \bigcup_{i=1}^\infty O_i $ 
differs from $O$ only by the set  of measure zero $\bigcup_{i=1}^\infty \partial U_i$. 
And since $O_i \subseteq U_i$ the first case gives 
$ \int_N \#({F_{O_i}}^{-1}(p)) \, d\mu_N(p) = \int_{O_i} |\Jac(F)| \, d\mu_M$. 
If we define $\Omega_i := \bigcup_{j=1}^i O_i$, then by the above observation and the disjointness 
of the $O_i$ we have:
$$ \eqalign{
 \int_N \#({F_{\Omega_i}}^{-1}(p)) \, d\mu_N(p)  &= \sum_{j=1}^i   \int_N \#({F_{O_j}}^{-1}(p)) \, d\mu_N(p) \cr
                                                                                   &= \sum_{j=1}^i  \int_{O_j} |\Jac(F)| \, d\mu_M \cr
                                                                                   &=  \int_{\Omega_i} |\Jac(F)| \, d\mu_M .
 }
 $$
 Now on the left, the sequence of integrands $\{\#({F_{\Omega_i}}^{-1}(p))\}$ converges monotonically  to 
 $\#({F_{O^\prime}}^{-1}(p))$, so by the Monotone Convergence (aka Beppo-Levi) Theorem 
 the limit of the sequence of integrals on the left is $\int_N \#({F_{O^\prime}}^{-1}(p)) \, d\mu_N(p)$,
which equals  $\int_N \#({F_O}^{-1}(p)) \, d\mu_N(p)$ since the two integrands agree almost everywhere.
 And since $O^\prime$  is the increasing union of the $\Omega_i$, the integrals on the 
 right converge to $\int_{O^\prime} |\Jac(F)| \, d\mu_M = \int_{O} |\Jac(F)| \, d\mu_M$, which completes 
 the proof that  $\int_N \#({F_O}^{-1}(p)) \, d\mu_N(p) = \int_{O} |\Jac(F)| \, d\mu_M$. \qed

\section{What Does it Mean to Distribute Points ``Uniformly''}

  Let $\Sigma$ denote a  compact, smoothly immersed hypersurface in $\Reals^n$, possibly with
smooth boundary.  To get started we need to make precise what it means to distribute 
points ``uniformly'' on $\Sigma$. Then, later, we will derive efficient algorithms for carrying this out. 
 Roughly speaking, to say that an algorithm selects points uniformly on a surface means that 
 the probability of selecting a point from a given subset should be proportional to the area of the 
 subset.  But giving a rigorous definition of this concept turns out to be considerably more involved than 
one might at first expect, and is of considerable independent interest, so we begin 
with a short discussion of the concept of equidistribution  in a more general context, 
with references for readers who would like to see more detailed accounts.  
We will frame this discussion in measure theoretic terms, however the reader who is unfamiliar 
with (or has forgotten) measure theory can, without serious loss of generality, interpret ``measure''
to mean the  usual volume measure on a smooth submanifold of $\Reals^n$.\par

\subsection{Equidistributed Sequences} 

\remark{Notation}
  In all that follows, $\mu$ will denote a finite, non-trivial Radon measure on a locally compact metric
space $X$ and $\mutilde$ will denote the corresponding probability measure ${1\over \mu(X)}\mu$. 
We will denote by $\{x_n\}$ a sequence of points  in $X$. For a subset $E$ of $X$ and positive integer 
$N$ we recall from section 2.2 that we denote by $\#(\{x_n\}_N, E)$ the number of integers $n$  with 
$1 \le n \le N$  such that $x_n \in E$.
\endremark
\smallskip

 A seemingly intuitive definition of the sequence $\{x_n\}$ being ``uniformly distributed'' in $X$ 
would be to demand that, for every measurable $E$, the ratio ${\#(\{x_n\}_N, E)\over N} $ 
should approach $\mutilde(E)$ as $N$ tends to infinity, or in the terminology of section 2.2, 
that $ \Prob(\{x_n\} \in E)$ is defined and equals $\mutilde(E)$.
But that condition could never be met by any 
sequence  $\{x_n\}$, as is easily seen by taking for $E$ the set of all the $x_n$.
That condition cannot hold even for all open sets $E$,
since we could take for $E$ the union of small open balls $B_n$ centered at the $x_n$ 
with $\mu(B_n) < \epsilon/2^n$ for arbitrarily small $\epsilon$. It turns out that we must restrict 
the subsets $E$ to those whose topological boundary
\footnote{i.e., the intersection of the closure of $E$ and the closure of its complement.}, 
$\partial E$,  has measure zero.

\definition{Definition 1} 
If $\mu$ is a finite and non-trivial Radon measure on a locally compact metric space $X$, then a 
sequence of points $\{x_n\}$ in $X$ is said to be {\it $1$-equidistributed\/} with respect to 
$\mu$  if it satisfies either and hence both of the following two equivalent conditions:

\item{1)} For every measurable set $E$  whose topological boundary has measure zero, \newline
 $\Prob(\{x_n\} \in E) $ is defined and equals ${\mu(E)\over \mu(X)} .$
 \smallskip

\item{2)}  
$${1 \over  \mu(X)}\int_X f(x) \, d \mu(x) = 
\lim_{N \to \infty} {1\over N} \sum_{n=1}^N f(x_n)$$
holds for every bounded continuous function $f: X \to \Reals$.
\enddefinition

\noindent   
\remark{Remark 1} If we take for the function $f$  the characteristic function of the set $E$
then, while of course such an $f$ will not in general be continuous, it is worth noting that in this case 
1) and 2) reduce to the same thing. 
\endremark 

\noindent   
\remark{Remark 2}  Whenever the equality 2) holds for a function $f$, we will say that the integral of 
the function $f$ can be evaluated by the Monte Carlo method, using the sequence $\{x_n\}$. 
\endremark 

\smallskip

\noindent   
\remark{Remark 3}  Recall that a sequence $\{\mu_n\}$ of finite Borel measures on $X$ is said
to {\it converge weakly\/} to the Borel measure $\mu$ if for every bounded continuous function $f: X \to \Reals$,
the sequence of integrals $\int f(x) \, d\mu_n(x)$ converges to $\int f(x) \, d\mu(x)$.  If $\delta_x$ denotes the 
Dirac delta measure at $x$ then 2) just says that the sequence of averages $ {1 \over N} (\delta_{x_1} + \ldots +\delta_{x_N})$ 
 converges weakly to $\mutilde$.

\endremark 

\smallskip

    The equivalence of the two conditions of Definition 1 is a consequence of a more general result of probability 
 theory referred to as the {\it Portmanteau Theorem\/}. We shall not give a proof here but refer instead 
 to Section 2, Chapter 1 of [B], where  Billingsley  proves the Portmanteau Theorem as his Theorem 2.1. 
 Other good sources for a careful discussion of this equivalence and related questions are [K1] and [K-N].
 In our case of primary interest, where $X$ is a smooth hypersurface embedded in some Euclidean space  
and $\mu$ is the area measure, we note that if the sequence $\{x_n\}$ is  $1$-equidistributed then in particular 
$\#(\{x_n\}_N, O)/N$ converges to $\mu(O)/\mu(X)$ for all open subsets $O$ of $X$ with smooth boundary, 
and this is a good intuitive way to think about the meaning of $1$-equidistributed in this setting.
 \smallskip

 It turns out that being $1$-equidistributed is a fairly weak restriction on the sequence $\{x_n\}$, and in particular it 
 has little to do with its being random. We can however use it to define a sequence of increasingly more 
 stringent conditions---called being $k$-equidistributed, and if  $\{x_n\}$ is $k$-equidistributed 
 for all positive integers $k$ then we will say that it is {\it completely equidistributed\/}. And being 
completely equidistributed {\it does\/} turn out to imply many of the conditions one would impose on a 
sequence for it to be called pseudo-randomly distributed with respect to $\mu$.  (For the case 
$X = [0,1)$ with $\mu$ Lebesgue  measure, see [F].) First let's look at the definition of $2$-equidistributed.
  Given an infinite sequence $\{x_1,x_2,\ldots\}$ 
in our measure space $(X,\mu)$ as above, we get an associated sequence 
$\{(x_1,x_2),(x_2,x_3), \ldots\}$ in the Cartesian product $X^2 :=X\times X$, and 
we say that the former sequence is $2$-equidistributed in $(X,\mu)$ if the latter sequence is 
$1$-equidistributed in $(X^2, \mu^2)$, where $\mu^2 := \mu \times \mu$ 
denotes the usual product measure. More generally:

\definition{Definition 2} 
If $X$, $\mu$, and $\{x_n\}$  are as in Definition 1 and $k \ge 2$ is a 
positive integer, then $\{x_n\}$ is said to be {\it $k$-equidistributed\/} with 
respect to  $\mu$  if the sequence in $X^k$ whose $n$-th
element is $(x_n,x_{n+1}\ldots, x_{n+k-1})$
is $1$-equidistributed with respect to $\mu^k$, i.e., 
if for each bounded, continuous real-valued function 
$f: X^k \to \Reals$ 
$${1 \over  \mu(X)^k }\int_X f(x_1.\ldots,x_k) \, d \mu(x_1)\ldots d \mu(x_k) = 
\lim_{N \to \infty} {1\over N} \sum_{n=1}^N f(x_n,x_{n+1}\ldots, x_{n+k-1}).$$
If $\{x_n\}$ in $X$ is {\it $k$-equidistributed\/} for all positive integers $k$ then 
we say that the sequence $\{x_n\}$ is {\it completely equidistributed\/} with 
respect to  $\mu$. 
\enddefinition

\noindent
So $\{x_n\}$ is $k$-equidistributed if and only if the Monte Carlo approach for evaluating 
$k$-fold integrals using the sequence $\{x_n\}$ works for all bounded continuous real-valued 
functions of $k$ variables on $X$, and the Portmanteau Theorem says that this is 
equivalent to being able to estimate the product measure of a subset $E$ of $X^k$ 
as the limit of the average number of points of the sequence 
$\{(x_n,x_{n+1}\ldots, x_{n+k-1})\}$ that belong to $E$, provided the boundary of $E$ 
has measure zero.

Since a function on $X^k$ can be regarded as a function on $X^{k+l}$ that is independent 
of its last $l$ arguments, it follows that a $k$-equidistributed sequence $\{x_n\}$ is also 
$k'$-equidistributed for all $k'$ with $1 \le k' < k$. However, a sequence that is 
$k$-equidistributed need {\it not\/} be $k'$-equidistributed for $k' > k$;  the following is an 
example of a $1$-equidistributed sequence that is not $2$-equidistributed.

\subsection{The van der Corput Sequences} 

  It is usually not easy to write down a simple and explicit equidistributed sequence. However 
there is an old (1935) and elegant family of examples for the unit interval $(0,1)$ that is due to 
J.~G.~van der Corput (for details, see [J] p.216).  Namely for a fixed positive integer greater than $1$ (the Radix) write the 
sequence of positive integers in that Radix, ending each with .0, and then ``reflect across the Radix point''. 
For example, if $\text{ Radix} = 2$, the sequence of integers in binary  notation is:
$$\text{\MySmallFont  1.0, 10.0, 11.0, 100.0, 101.0, 110.0, 111.0, 1000.0, 1001.0, 1010.0, 
1011.0, 1100.0, 1101.0, 1110.0, 1111.0, }\ldots$$ 
and reflecting each element of this sequence in its binary point gives:
$$\text{\MySmallFont  0.1, 0.01, 0.11, 0.001, 0.101, 0.011, 0.111, 0.0001, 0.1001, 0.0101, 
0.1101, 0.0011, 0.1011, 0.0111, 0.1111, }\ldots$$ 
or using fractional notation:
$$ \tfrac{1}{2}, \ \ \ \;  \tfrac{1}{4}, \ \ \ \; \tfrac{3}{4}, \ \ \ \; \tfrac{1}{8}, \ \ \  \;\tfrac{5}{8}, \ \ \ \; 
\tfrac{3}{8}, \ \ \ \; \tfrac{7}{8}, \ \ \ \;
\tfrac{1}{16}, \ \ \ \; \tfrac{9}{16}, \ \ \ \; \tfrac{5}{16}, \ \ \ \;\tfrac{13}{16}, \ \ \ \; \tfrac{3}{16}, \ \ \ \; \tfrac{11}{16}, 
\ \ \ \; \tfrac{7}{16}, \ \ \ \; \tfrac{15}{16}, \ \ \  \; \ldots $$
and this is, by definition,  the Radix = 2 or binary van der Corput sequence.
This recipe for forming the van der Corput sequences looks a bit mysterious, and, naively,  it might appear 
more natural to rearrange the above sequence in the following order:
$$ \tfrac{1}{2}, \ \ \ \;  \tfrac{1}{4}, \ \ \ \; \tfrac{3}{4}, \ \ \ \; \tfrac{1}{8}, \ \ \  \;\tfrac{3}{8}, \ \ \ \; 
\tfrac{5}{8}, \ \ \ \; \tfrac{7}{8}, \ \ \ \;
\tfrac{1}{16}, \ \ \ \; \tfrac{3}{16}, \ \ \ \; \tfrac{5}{16}, \ \ \ \;\tfrac{7}{16}, \ \ \ \; \tfrac{9}{16}, \ \ \ \; \tfrac{11}{16}, 
\ \ \ \; \tfrac{13}{16}, \ \ \ \; \tfrac{15}{16}, \ \ \  \; \ldots $$
i.e., arrange each successive sub-series of proper fractions with denominator $2^n$ and odd numerators in increasing 
order.  However, while the original binary van der Corput sequence is particularly well equidistrbuted 
(technically, it is what is called a ``low-discrepancy sequence'') the rearranged sequence turns out {\it not\/} to be 
equidistributed at all!

\section{Random Number Generators (RNGs)}

   The further theory of equidistributed sequences is based on so-called {\it Random Number 
Generators\/} (or RNGs), which we shall always understand here to mean highly equidistributed 
sequences $\{\xi_j\}$ in the interval $[0,1)$ of  $\Reals$ or, more precisely, algorithms for constructing 
such sequences. From its inception, this subject has drawn the attention of many famous mathematicians, 
including Herman Weyl [W],  Stanislaw Ulam [U], John von Neumann [VN], and Donald Knuth [K1]. There is an 
extensive literature devoted to the subject, and many excellent articles and books that provide careful 
and detailed coverage of the theory, including:  [B], [B-C], [E], [K1], [K-N]. Here we will only describe the basic 
facts concerning RNGs, concentrating on those we will need later in this article, and for proofs and 
further details we will refer to the above sources.

    An obvious first question to consider is whether completely equidistributed RNGs even exist and if so, 
how to go about creating them. If one is satisfied with the pure mathematician's usual non-constructive approach to 
the real numbers the answers to these questions seem encouraging. It was shown already in [F] that for almost all real 
$\theta$ larger than one, the fractional part of $\theta^n$ defines a completely equidistributed RNG, and while 
no explicit such $\theta$ has been identified, in [K1], Knuth does describe an explicit completely equidistributed RNG.

  But RNGs play an essential role in practical computer science, making it important to have RNG algorithms 
 that can be coded efficiently on modern computers. Here the situation is considerably more complex.  
 As is well-known, naive attempts at a rigorous machine-oriented theory of RNGs quickly 
 founder on apparent paradoxes. Since even a $1$-equidistributed sequence is dense, it must {\it a fortiori\/} 
 be infinite, and in particular it cannot be periodic. But if one tries to design an algorithm to produce a RNG 
 on a computer using a fixed precision to represent real numbers, there are only a finite set of reals that can 
 be represented. And even if one uses ``arbitrary precision arithmetic'' (i.e., allows the number of computer words 
 used to represent a real number to grow as more precision is required) then although there are potentially infinitely 
 many real numbers available, the sequence will still be periodic. This is because any real machine has only finitely 
 many states, and once the algorithm has produced enough elements of the sequence, it must eventually repeat a 
 state of the machine and at this point the sequence will become periodic.  There is no panacea for solving this 
 problem, but a careful discussion of how to approach it will be found in the first section of L'Ecuyer's  review 
 article [E]. Suffice it to say that for any realistic application it is possible to find a very highly equidistributed 
 RNG that works ``for all practical purposes''.  To give some idea of what is available in this direction, there is 
 the remarkable RNG, called The Mersenne Twister, due to Matsumoto and Nishimura [M-N], that has a period of 
 $2^{19937}-1$ and, to 32-bit accuracy, is $623$-equidistributed!

\section{Methods for Constructing New Equidistributed Sequences from Old}

  For the remainder of this article we assume that we have available a $k$-equidistributed RNG 
for some large $k$. We will denote it by $\{\xi_n\} = \{\xi_1, \xi_2,\ldots\}$, and use it to help devise 
methods for constructing highly equidistributed sequences in various spaces. Here are some basic examples.

\remark{Method 1. Products of Intervals}\   If $a_1,\ldots,a_n$ and $b_1,\ldots, b_n$ are two $n$-tuples 
of real numbers with $a_i < b_i$, then we can construct a $(k - n)$-equidistributed sequence in the  product 
$[a_1,b_1)\times \ldots \times [a_n,b_n)$ by taking the $j$-th element to be
$(a_1 + (b_1 -a_1)\xi_j,\ldots, a_n + (b_n -a_n)\xi_{j+n-1} )$. This follows trivially from the definition of 
$k$-equidistribted.
\endremark

\remark{Method 2. Subsets; The Method of Rejection}\  Suppose that in our measure space $(X,\mu)$ 
we have a $k$-distributed sequence $\{x_i\}$ and let $Y$ be a measurable subset of $X$ having positive 
measure and whose boundary $\partial Y$ has measure zero. Define a sequence $\{y_i\}$ in $Y$ by taking 
$y_j$ to be the $j$-th element of the sequence  $\{x_i\}$ that belongs to $Y$;  i.e., simply ``reject'' all elements 
of the original sequence that do not lie in $Y$. Then it is immediate from Definition 1 that the sequence 
$\{y_i\}$ is $k$-equidistributed in $Y$.
\endremark

\remark{Method 3. Unions}\  We can also go the other way. That is, suppose we have a $k$-equidistributed 
sequences $\{ x^i_j\}$ for each of $N$ disjoint measurable subsets $E^i$ of $(X,\mu)$. We can 
construct a $k$-equidistributed sequence $\{x_j\}$ for their union $E$ as follows. Define  
$M_n := \sum_{i=1}^n \mu(E_i)$. To obtain $x_j$,  first choose one of the $E_\ell$ 
with a probability equal to its relative measure by seeing which subinterval $[M_\ell,M_{\ell+1}]$ of 
$[0,M_N)$ contains $M_n \xi_j$, and let $x_j$ be the first element of the sequence $\{x^\ell_i\}$ that 
has not already been chosen.
\endremark

\remark{Method 4. Measure Preserving Maps}\   Assume that we have a $k$-equidistributed sequence $\{ x_n\}$ 
for $(X,\mu)$, a second measure space $(Y,\nu)$, and a continuous map $f$ of $X$ onto $Y$ such 
that  $\mu(f^{-1}(E))$ is proportional to $\nu(E)$. Then the sequence $\{ f(x_n)\}$ is clearly $k$-equidistributed in $Y$.
\endremark

\remark{Method 5. Equidistributed Sequences for the ball, $\B^n$ and Sphere $\S^{n-1}$}\  
We next show how to  combine the above methods to create $(k-n)$-equidistributed 
sequence $\{b_j\}$ for the $n$-dimensional ball $\B^n$ and $\{s_j\}$ for the $(n-1)$-sphere $S^{n-1}$.  
We start by using Method 1 to define a $(k-n)$-equidistributed sequence $\{x_j\}$ for 
the cube $C_n := [-1,1]^n$; namely $x_j := ( 2 \xi_j - 1, \ldots, 2 \xi_{j+n-1} -1 )$. The (punctured) 
unit ball $\B^n := \{ x \in \Reals^n \mid 0 <||x|| < 1\}$ is a subset of $C_n$ with smooth boundary, 
so we can apply rejection (Method 2)  to the sequence $\{x_j\}$ to get a $(k-n)$-equidistributed 
sequence $\{ b_j\}$ for $\B^n$. The normalizing map $f: \B^n \to \S^{n-1}$, defined by $f(x) := {x \over ||x||}$, 
satisfies the condition of Method 4, that is, if $E$ is a measurable subset of $\S^{n-1}$, its surface 
measure is proportional to the volume measure of $f^{-1}(E)$, the cone over $E$, hence 
$s_j := {b_j \over ||b_j||}$ is a $(k-n)$-equidistributed sequence for $\S^{n-1}$.

{\bf A Small Problem!}  A careful analysis exposes a serious difficulty with the simple and elegant
Method 5  when $n$ becomes large. Note that the probability of the point $x_j \in C_n$ being 
``accepted'' (i.e., not rejected ) as an element of $\B^n $ is the ratio $\rho_n :=\kappa_n/2^n$ of the 
volume $\kappa_n$ of the $n$-ball to the volume $2^n$ of $C_n$.  When $n=3$ this ratio is 
$\rho_3=(4\pi/3)/8 = \pi/6$, so slightly more than half of the $x_j$ will be accepted. But this is deceptive 
as a clue to what happens for $n$ large. There is a well-known formula for $\kappa_n$, namely 
$\kappa_n = \pi^{n/2}/ \Gamma({n\over 2}+1)$. Restricting for simplicity to the case of even $n = 2k$,  
we have $\kappa_{2k} = \pi^k /k!$ so $\rho_{2k} = (\pi/4)^k/k!$, and we see that, already 
for $n = 10$, only about one in four hundred of the $x_j$ will be accepted! Moreover the situation gets rapidly 
worse as $n$ increases, making Method 5 impossibly slow, so we need another, faster way
to select random points of $\B^n$ and $\S^{n-1}$ for large $n$. We will consider only the case of $\S^{n-1}$, 
since it is easy to see that if $x$ is a random point of $\S^{n-1}$ and we choose a random $u \in [0,1)$ 
with our RNG, then $x \cdot u^{1/n}$ is a random point of $\B^n$.

  We recall that the Normal (or Gaussian) Distribution with mean $\mu$ and variance $\sigma^2$ is the 
continuous probability distribution (or density) on $\Reals$ given by 
$$
f(x) := {1\over {\sigma \sqrt{2 \pi}}}e^{-{(x-\mu)^2\over 2 \sigma^2}}
$$
and we will denote by $\Cal N(\mu,\sigma)$ the associated measured manifold. In particular, 
$\Cal N(0,1)$ is often referred to as the Standard Normal Deviate. There are numerous methods for 
starting from a highly equidistributed RNG  and using it to construct a sequence of Standard Normal 
Deviates, that is a sequence $\{x_i\}$ of positive real numbers that are independent and are distributed 
according to the distribution $\Cal N(0,1)$; see for example Algorithm M of [K1] page 110. 
Then it is easy to check that the sequence $y_k := (x_k, x_{k+1},\ldots,x_{k+n})$ of points of $\Reals^n$ 
has a distribution that is invariant under rotation, and it follows that the normalized sequence 
$\{y_k/||y_k||\}$ is highly equidistributed in $\S^{n-1}$. For details, see [K1] page 116.

\endremark

\section{The Space of Lines and its Kinematic Measure}

\subsection{The space $\L^n$ of Oriented Lines in $\Reals^n$}\  
  We will denote by $\L^n$ the space of oriented lines in $\Reals^n$.
A parametric equation for some $\ell \in \L^n$ has the form $x(t) =  t v+p$, 
where $p$ is any point of $\ell$ and $v$ is a non-zero vector giving its direction. 
Replacing $p$ by any other point on $\ell$ or multiplying $v$ by a positive 
scalar gives another parametric equation for the same oriented line,  
so the mapping from $\L^n$ to pairs $(v,p)$ is highly non-unique, and since we wish
to identify each $\ell$ with a specific pair $(v,p)$ we next specify how to choose 
$v$ and $p$ canonically. We make $v$ unique by normalizing it, i.e., we 
replace $v$ by ${v\over ||v||}\in \S^{n-1}$. This is equivalent to 
making the parameter $t$ the signed arc-length measured from $p$. And to make 
$p$ unique we choose it to be the point of $\ell$ nearest to the origin, which is equivalent 
to requiring $p$ to be orthogonal to $v$. 
In this way we identify $\L^n$ with the set of all pairs $(v,p)$ with $v\in \S^{n-1}$ 
and $p \in v^\perp$, the $n-1$ dimensional subspace of $\Reals^n$ orthogonal to $v$.  
\footnote{Note that $v^\perp$ is  just $T_v \S^{n-1}$, the tangent space to $\S^{n-1}$ at $v$, so we 
can identify $\L^n$ with  $T \S^{n-1}$, the tangent bundle of $\S^{n-1}$, however 
we shall not have occasion to make use of this identification.}
To each such $(v,p)$ we associate the oriented parametric line $\sigma^{(v,p)}(t) := t v +p$, 
and we note that $p$ is the orthogonal projection of each point of this line onto $v^\perp$.
We will denote by $ \L^n_r$ the subset of $ \L^n$ consisting of lines that meet 
the ball $\B^n_r$ of radius $r$ centered at the origin. Clearly $\ell = (v,p)$ belongs to 
$ \L^n_r$ if and only if $||p|| < r$.

 \subsection{The Action of Euclidean Motions on $\L^n$}\  
The group of Euclidean motions of $\Reals^n$ (transformations that
can be written as a rotation followed by a translation) acts in an obvious way on
oriented lines, and this action is transitive, since clearly any oriented line 
can  be mapped onto any other by a rotation followed by a translation. 
A rotation $g$ maps $\sigma^{(v,p)}$ to $\sigma^{(gv,gp)}$, so
for $(v,p) \in \L^n$,  $g(v,p) = (gv,gp)$.
But how does a translation $\tau_a: x \mapsto x+a$ act on the line $(v,p)$? 
Clearly translation does not affect the direction component, $v$, and translation 
by a multiple of $v$ does not affect $p$ either, 
so if we write $a_{v^\perp}$ for the component of $a$ 
orthogonal to $v$, then $\tau_a(v,p) = (v, p + a_{v^\perp})$. 
In particular, on $\{v\} \times v^\perp$ (the space of all lines with direction $v$)
$\tau_a$ is just translation by $a_{v^\perp}$ on the second component, 
and so in particular it is Lebesgue measure preserving.
 
 \subsection{The Kinematic Measure for $\L^n$} \ 
  We next note that there is a natural measure on $\L^n$ that we denote
by $\mu$ and call {\it kinematic measure\/}. To specify it, we can use the  
identification of Radon measures with continuous linear functions on the Banach 
space of bounded, continuous functions with compact support. That is, given any 
bounded continuous function with compact support $f(v,p)$ on $\L^n$,  
we must define its integral, $\int f(v,p) \, d\mu$. 
For each $v \in \S^{n-1}$, first integrate $f$ over $v^\perp$ with respect to Lebesgue
measure, getting in this way a function $F$ on  $\S^{n-1}$; $F(v) := \int_{v^\perp} f(v,p) \, dp$,  
and then integrate $F(v)$ over $\S^{n-1}$ with respect to its usual area measure
$d \sigma$, so the final result is 
$$\int_{\L^n}f(v,p) \, d\mu := \int _{\S^{n-1}} d\sigma(v)\int_{v^\perp} f(v,p)\, dp.$$
We claim that both translations and rotations, and so the whole group of Euclidean 
motions, preserves  kinematic measure. In fact,  we saw above that 
translations act on $\L^n$ by a translation in each $v^\perp$, and 
so they clearly preserve kinematic measure. And we saw that rotations act by $g(v,p) = (gv,gp)$, 
from which it follows they also preserve the kinematic measure, since $g$ maps $v^\perp$ 
isometrically onto $(gv)^\perp$ and is also an isometry of $\S^{n-1}$.  It is 
well-known (and elementary) that if a group acts smoothly and transitively on a 
smooth manifold $M$, then a non-trivial invariant volume form for $M$, if one exists,
is unique up multiplication by a positive scalar. So in particular kinematic measure is 
characterized up to a positive scalar multiple by the fact that it is invariant under 
Euclidean motions.

\subsection{Constructing A Highly Equidistributed Sequences of Lines} 

  Starting  from our $k$-equidistributed RNG, $\{\xi_n\}$, we can now construct a sequence  
$\{\ell_j = (v_j,p_j) \}$  in  $\L_r^n$ that is $(k-2n)$-equidistributed with respect to 
kinematic measure. In fact, it is clear from the definition of kinematic measure that
we will get such a sequence if we use Method 5 of Section 5 to first select $v_j \in \S^{n-1}$ 
and then to also select $p_j$ in the ball of radius $r$ of $ v_j^\perp$.  \par

   Below we will explain the final steps that converts the sequence $\{\ell_j\}$ into a $(k-2n)$-equidistributed 
sequence of points $\{x_j\}$ on any hypersurface $\Sigma \subset \B^n_r$, essentially  by intersecting 
the lines $\ell_j$ with $\Sigma$. But for a rigorous justification of this approach we must first recall 
some basic facts about transversality and The Cauchy-Crofton Formula.

\section{Transversality}                     

  As just mentioned, a crucial step in our technique for constructing an equidistributed sequence of points 
on a bounded hypersurface $\Sigma \subset \Reals^n$  involves producing finite subsets of $\Sigma$ 
by intersecting it with a line $\ell$. Intuitively we expect that a line should intersect a hypersurface in a 
discrete subset, and so a finite one if the hypersurface is bounded.
But we must be careful; $\ell \cap \Sigma$ {\it can\/} be infinite, for example if $\Sigma$ is part of a 
cone or cylinder $C$ and $\ell$ is a generator of $C$. To make our technique rigorous we make a  
short excursion into transversality theory, where we will find that $\ell \cap \Sigma$ is indeed finite 
except possibly for a set of lines $\ell$ having kinematic measure zero. First we recall:

\definition{Definition} Let $M$, $N$ be smooth manifolds and $\Sigma$ a smooth submanifold of $N$.
We call a smooth map $f: M \to N$  {\it transverse to\/} $\Sigma$ if for all $p \in f^{-1}(\Sigma)$, 
$\Im(Df_p)$ together with $T\Sigma_{f(p)}$  span $TN_{f(p)}$.
\enddefinition
 
\remark{Remark}   We note that if $f$ is a submersion, that is if for each $p \in M$, $Df_p$ maps 
$TM_p$ onto $TN_{f(p)}$, then $f$ is automatically transverse to any submanifold of $N$.
\endremark

\proclaim{Theorem} If a smooth map $f: M \to N$ is transverse to a submanifold $\Sigma$ of $N$ then 
$f^{-1}(\Sigma)$ is a smooth submanifold of $M$ whose codimension in $M$ is the same as the codimension 
of $\Sigma$ in $N$.
\endproclaim
\proof This is an elementary result; see Theorem 3.3 on page 22 of [H].  \qed
\remark{Example}
If $\ell = (v,p) \in \L^n$ and the map $\sigma^{(v,p)} : \Reals \to \Reals^n$, $t \mapsto t v + p$ is 
transverse to a bounded hypersurface $\Sigma$ of $\Reals^n$, then the theorem says that the set of 
$t \in \Reals$ such that $\sigma^{(v,p)}(t) $ belongs to $\Sigma$ is a zero dimensional submanifold of $\Reals$,
i.e., a discrete subset of $\Reals$, hence its image, $\ell \cap \Sigma$, is a discrete and hence finite 
subset of $\Sigma$.
\endremark

  We recall that a subset of a smooth manifold is called a {\it null set\/} if it meets each chart 
neighborhood in a subset that has Lebesgue measure zero in that chart, or equivalently if it has 
measure zero with respect to a measure defined by a non-vanishing differential form of top degree. 
A property is said to hold for almost all points if it holds outside of a null set.

\proclaim{Parametric Transversality Theorem} Let $M$, $N$ and $P$ be smooth manifolds, 
$F: M \times P \to N$ a smooth map and $\Sigma$ a smooth submanifold of $N$. Assume 
that $F$ is a submersion, or more generally that it is transverse to $\Sigma$. If for $p \in P$ 
we denote by $F_p : M \to N$ the map $m \mapsto F(m,p)$, then $F_p$ is transversal to 
$\Sigma$ for almost all $p$ in $P$.
\endproclaim
\proof This is Theorem 2.7 page 79 of [H].  \qed

  We next apply this with $M = \Reals$, $N = \Reals^n$, $P =\L^n$ 
as the parameter space, and defining $F: \Reals \times \L^n \to \Reals^n$ by 
$F(t,(v,p)) := \sigma^{(v,p)}(t) = tv +p$. If $q = (t,(v,p)) \in \Reals \times \L^n$, recall that 
$v \in \S^{n-1}$ and $p \in v^\perp$, so 
then $T( \Reals \times \L^n)_q = T(\Reals)_t \oplus T(\S^{n-1})_v \oplus v^\perp$, and 
 $DF_q({\partial \over \partial t}) = v$ while for $w \in v^\perp$, $DF_q(w) = w$, so the image 
 of $DF_q$ includes all of $\Reals v \oplus v^\perp = \Reals^n$, proving that $F$ is a submersion. 
 Hence the above example together with the Parametric Transversality Theorem imply that:
 
 \proclaim{Corollary} If $\Sigma$ is a smooth hypersurface in $\Reals^n$, then all $\ell \in \L^n$ 
 except for a set of kinematic measure zero intersect $\Sigma$ in a discrete set, and hence in a
 finite set if $\Sigma$ is bounded. 
 \endproclaim
 
 \section{Classic Cauchy-Crofton}
 
  Many interesting geometric quantities associated to submanifolds of $\Reals^n$ can be expressed 
as an integral of some related function over $\L^n$ with respect to kinematic measure. The study of
such relations is known as {\it Integral Geometry\/}, and one of its earliest results is a formula of Crofton 
giving the length of a plane curve in terms of the integral over $\L^2$ of the number of times a line 
$\ell$ intersects it.  Both Crofton's result and its proof generalize to a formula for the area of a compact 
hypersurface of $\Reals^n$---and this generalization,  which we will call The Classic Cauchy-Crofton 
Formula, plays an important role in our algorithm for constructing point clouds on hypersurfaces.

  In what follows $\Sigma$ denotes a compact, smooth hypersurface in $\Reals^n$. For $\ell \in \L^n$ 
we define $\#(\ell \cap\Sigma)$  to be zero if $\ell$ does not meet $\Sigma$ or if it meets it non-transversally, 
and otherwise it is the number of points in $\ell \cap\Sigma$.  We recall that we write the  volume $(n-1)$-form 
of a hypersurface $\Sigma$ as $d\omega_\Sigma$ and  its integral, $A(\Sigma)$, over $\Sigma$ is called the area 
of $\Sigma$. Finally, we will denote by $\kappa_m$ the volume of the $m$-ball $\B^m$.
\footnote{We recall the well-known formula, $\kappa_m = {\pi^{m\over 2}   /   \Gamma\left({m\over 2} + 1\right)}$.}
\smallskip

\proclaim{The Classic Cauchy-Crofton Formula} 
\footnote{Most integral geometry sources omit the $2$ in the denominator; this is because they use the 
space of unoriented lines rather than the space of oriented lines as we do here.}
  If $\Sigma$ is a compact, smooth hypersurface in $\Reals^n$, then 
$$A(\Sigma) = {1\over 2 \kappa_{n-1}} \int_{\L^n}  \#(\ell \cap\Sigma)\, d\mu(\ell) .$$
 \endproclaim
 
 A standard reference for the proof of Crofton's original formula, for a curve in $\Reals^2$,
 is  pages 12--13 of Santalo's classic {\it Introduction to Integral Geometry\/} [S1]. \par

  We begin our proof of Classic Cauchy-Crofton with a geometrically intuitive observation: 
If a cylindrical  flashlight beam meets a flat surface at an angle $\theta$, then the area $A(S)$ of the region 
$S$ of the surface illuminated by the beam is related to the area $A(\tilde S)$ of the beam cross-section $\tilde S$ 
by $A(\tilde S) = A(S)| \cos(\theta)|$. The following is a generalization of this fact to higher 
dimensions, using more precise language.  \par 

\proclaim {Proposition 1}
  Let $\nu, v \in \S^{n-1}$ and let  $\Pi_v$ denote orthogonal projection of $\Reals^n$ onto $v^\perp$,  so that  
$\Pi_v(x) := x - \ip(x,v)v$. If $\pi_v: \nu^\perp \to v^\perp$ denotes the restriction of $\Pi_v$ to $\nu^\perp$, 
then the Jacobian of $\pi_v$ has the constant value $ |\ip(v,\nu)| = |\cos(\theta_v)|$, 
where $\theta_v$ is the angle between $\nu$ and $v$. 
(Equivalently,  If $S$ is a bounded region of $\nu^\perp$, then the area $A(S)$ of $S$ and the area  
$A(P_v(S) )$ of its projection onto $v^\perp$ are related by:
$A(P_v(S)) = A(S) |\ip(v,\nu)| = A(S) |\cos(\theta_v)| )$.
\endproclaim

\proof
The case $v = \nu$ is trivial, so we assume  $v \ne \nu$,  identify $\Reals^2$ with the linear 
span of $\nu$ and $v$ and $\Reals^{n-2}$ with its orthogonal complement, $\nu^\perp \cap v^\perp$. 
Then $\Reals^n$ is the orthogonal direct sum $\Reals^2 \oplus \Reals^{n-2}$ and $\Pi_v$ is the identity 
on $\Reals^{n-2}$. This effectively reduces the Proposition to the case $n = 2$ where it is obvious.  \qed \par

   If we regard $\nu$ as the north pole, then the equatorial plane $\nu^\perp$ splits the unit sphere 
$\S^{n-1}$ into two hemispheres $\S^{n-1}_{\pm}$, and each of  $\phi_\pm := P_\nu | \S^{n-1}_{\pm}$ 
is a diffeomorphism of $\S^{n-1}_{\pm}$ onto the unit ball $\B^{n-1}$ in that plane. Since the tangent
space to $\S^{n-1}$ at a point $v$ is $v^\perp$, it follows from the calculation above that the absolute value 
of the Jacobian determinant of $\phi_\pm$ at any point $v$ is just  $|\ip(v,\nu)|$. By the change of variable 
formula of multi-variable calculus, this proves the following:

\proclaim{Proposition 2}
For all $\nu \in \S^{n-1}$ 
$$ \int_{\S^{n-1}} |\ip(v,\nu)| \, d\sigma(v) = 2 \kappa_{n-1},$$ where as in 
section 3, $d \sigma$ denotes the Riemannian volume element of $\S^{n-1}$.
\endproclaim

 We are now ready to complete the proof of the Classic Cauchy-Crofton Formula. 
 As above, for $v\in \S^{n-1}$, we denote orthogonal projection of $\Reals^n$ onto $v^\perp$ by $\Pi_v$, 
 and we will denote by $P_v: \Sigma \to v^\perp$ its restriction to $\Sigma$. For $x \in \Sigma$ we 
 denote by $\nu_x \in \S^{n-1}$ the unit normal to $\Sigma$ at $x$, so ${\nu}_x^\perp$ is the tangent 
 space to $\Sigma$ at $x$.
 
 \proclaim {Lemma} 
  $|\Jac(P_v)(x) |= |\ip(v,\nu_x)|$.  
 \endproclaim
 
 \proof
 Since $\Pi_v$ is linear, its restriction $\pi: {\nu}_x^\perp \to v^\perp$ to  the tangent space 
$ {\nu}_x^\perp$ of $\Sigma$ at $x$ is the differential $DP_x$ of the restriction $P$ of $\Pi_v$ to $\Sigma$,
so the Lemma follows from Proposition 1.
 \qed
 
  Then, recalling from Section 6 the definition of the kinematic measure on the space $\L^n$ of lines, and noting 
  that if $\ell = (v,p) \in \L^n$ then $(\ell \cap \Sigma) = {P_v}^{-1}(p)$,

   $$\eqalign{ 
    \int_{\L^n}  \#(\ell \cap\Sigma)\, d\mu(\ell) 
  &= \int _{\S^{n-1}} d\sigma(v)\int_{v^\perp} \#({P_v}^{-1} (p))\, dp    \cr
 &=  \int _{\S^{n-1}} d\sigma(v) \int_\Sigma|\Jac(P_v)|(x) \, d\Sigma(x)  \hbox{\hskip 20 pt by the Co-Area Formula of Section 2} \cr
 &= \int _{\S^{n-1}} d\sigma(v) \int_\Sigma |\ip(v ,\nu_x)|\, d\Sigma(x)  \hbox{\hskip 35 pt by the Lemma}   \cr 
 &= \int_\Sigma  \, d\Sigma(x) \int_{\S^{n-1}}  |\ip(v ,\nu_x)|\,  d\sigma(v)  \hbox{\hskip 35 pt by Fubini}   \cr
  &= A(\Sigma) \times  (2 \kappa_{n-1})  \hbox{\hskip 95pt  by Proposition 2} \cr
 }$$
 which completes he proof of the Classic Cauchy-Crofton Formula. 
\qed

 \section {The Cauchy-Crofton Principle}  
   We will need a generalization of The Classic Cauchy-Crofton Formula that we will call The Cauchy-Crofton Principle.
It says, roughly, that to integrate a function over a hypersurface $\Sigma$ of $\Reals^n$ you can first sum its values at
 the points where a line $\ell$ meets $\Sigma$, and then integrate that sum over the space $\L^n$ of  all the lines
 in $\Reals^n$. When $f \equiv 1$ this of course is just a restatement of  Classic Cauchy-Crofton and, as we will
 see in the next section, conversely The Cauchy-Crofton Principle  follows from Classic Cauchy-Crofton.

 \definition{Definition}
 If $f:\Sigma \to \Reals$ is a real-valued function on a compact hypersurface $\Sigma$ in $\Reals^n$
 we define $f^\Sigma : \L^n \to \Reals$ by 
 $f^\Sigma(\ell) := 0$  if $\ell$ is not transverse to $\Sigma$, and otherwise 
 $f^\Sigma(\ell) := \sum_{p \in \ell \cap \Sigma} f(p)$. 
 More generally, if $f: \Sigma^k\to \Reals$, we define $f^\Sigma : {(\L^n)^k} \to \Reals$ by 
 $f^\Sigma (\ell_1, \ldots, \ell_k) := 0$ if some $\ell_i$ is not transverse to $\Sigma$, and otherwise 
 $$f^\Sigma (\ell_1, \ldots, \ell_k) := \sum f (x^1_{j_1},\ldots, x^k_{j_k}),$$  
 where the sum is over all $k$-tuples $(x^1_{j_1},\ldots, x^k_{j_k})$ such that 
 $x^i_{j_i} \in \ell_i \cap \Sigma$.
 \enddefinition

\proclaim {The Cauchy-Crofton Principle}  
   If $\Sigma$ is a compact, smooth hypersurface in $\Reals^n$, and $f: \Sigma \to \Reals$ is a continuous 
real-valued function on $\Sigma$,  then
$$\int_\Sigma f(x) \,  d\omega_\Sigma(x) = {1\over 2 \kappa_{n-1}} \int_{\L^n} f^\Sigma(\ell)\, d\mu(\ell).$$
\endproclaim

\proclaim {Corollary: The Multivariable Cauchy-Crofton Principle}  
   If $\Sigma$ is a compact, smooth hypersurface in $\Reals^n$ and if $f: \Sigma^k \to \Reals$ is a continuous 
real-valued function on $\Sigma^k$,  then
$$\int_{\Sigma^k} f(x_1,\ldots, x_k) \,  d\omega_\Sigma(x_1)\ldots d\omega_\Sigma(x_k) = 
{\left(1\over 2 \kappa_{n-1}\right)}^k  \int_{{(\L^n)}^k} f^\Sigma(\ell_1, \ldots, \ell_k) \, d\mu(\ell_1)\ldots d\mu(\ell_k).$$
\endproclaim
\noindent
{\it Proof of the Corollary.\/}
We give the details for the case $k = 2$; the general case is only notationally more complicated. 
If for $x \in \Sigma$ we define $f_x:\Sigma \to \Reals$ by $f_x(x_2) := f(x, x_2)$, then  applying 
The Cauchy-Crofton Principle with $f = f_{x_1}$,  
$\int _\Sigma f_{x_1}(x_2) \, d\omega_\Sigma(x_2) = 
{1\over 2 \kappa_{n-1}} \int_{\L^n} f_{x_1}^\Sigma(\ell_2) \, d\mu(\ell_2) \hbox{,\ so} $
$$\eqalign{
\int_{\Sigma^2} f(x_1,x_2) \,  d\omega_\Sigma(x_1)\,d\omega_\Sigma(x_2) = 
\int_{\Sigma}  \, d\omega_\Sigma(x_1) \int_\Sigma f_{x_1}(x_2) \, d\omega_\Sigma\,d\omega_\Sigma(x_2)\cr
={1\over 2 \kappa_{n-1}}\int_{\Sigma}  \, dA(x_1) \int_{\L^n} f_{x_1}^\Sigma(\ell_2) \, d\mu(\ell_2),
}$$

or by Fubini's Theorem, 
$$\int_{\Sigma^2} f(x_1,x_2) \,  dA(x_1)\,dA(x_2) = 
{1\over 2 \kappa_{n-1}}\int_{\L^n}  \, d\mu(\ell_2) \int_{\Sigma}  f_{x_1}^\Sigma(\ell_2) \, dA(x_1) , $$
and applying The Cauchy-Crofton Principle again, this time with $f = g_{\ell_2}$ where 
$g_{\ell_2}(x_1) :=  f_{x_1}^\Sigma(\ell_2), $
$$\int_{\Sigma^2} f(x_1,x_2) \,  dA(x_1)\,dA(x_2) =  
{\left(1\over 2 \kappa_{n-1}\right)}^2\int _{\L^n}\, d\mu(\ell_2)  \int_{\L^n} g_{\ell_2}^\Sigma(\ell_1) \, d\mu(\ell_1), $$
and appealing to Fubini's Theorem again we only have to check that $ g_{\ell_2}^\Sigma(\ell_1)  = f^\Sigma(\ell_1,\ell_2)$. 
But 
$$g_{\ell_2}(x) = f_x^\Sigma(\ell_2) = \sum_{q \in (\ell_2 \cap \Sigma)} f_x(q) =
\sum_{q \in (\ell_2 \cap \Sigma)} f(x, q)  \hbox{\  so,} $$ 
$$g^\Sigma_{\ell_2}(\ell_1) = \sum_{p \in (\ell_1 \cap \Sigma)} g_{\ell_2}(p) =
\sum_{p \in (\ell_1 \cap \Sigma)} \sum_{q\in (\ell_2 \cap \Sigma)}f(p,q)  = f^\Sigma(\ell_1,\ell_2). \qed$$

\section{Demonstration of The Cauchy-Crofton Principle} 

  Our proof of The Cauchy-Crofton Principle from The Classic Cauchy-Crofton Formula will proceed by 
a series of reductions.  We note first that since the function $f^\Sigma: \L^n \to \Reals$ in its statement  clearly depends 
linearly on the function $f:\Sigma \to \Reals$ it follows that:

\proclaim {Linearity Lemma} If  The Cauchy-Crofton Principle holds for functions $f_\alpha:\Sigma \to \Reals$ then
it also holds for any finite linear combinations of the $f_\alpha$.
\endproclaim

 \proclaim {Proposition 1} If a compact, smooth hypersurface $\Sigma$ of $\Reals^n$ is a finite union of the interiors of  
 hypersurfaces $\Sigma_i$ and if each $\Sigma_i$ satisfies The Cauchy-Crofton Principle, then so does $\Sigma$. 
\endproclaim

\proof
  Let $\phi_i : \Sigma \to \Reals$ be a smooth partition of unity for $\Sigma$ with the support of $\phi_i$ a 
compact subset of the interior of $\Sigma_i$. Given a continuous $f: \Sigma \to \Reals$, the functions $f \phi_i$ 
satisfy The Cauchy-Crofton Principle and by the Linearity Lemma so does their sum $f$.
\qed
  
  This shows that the Cauchy-Crofton Principle is local in nature;  i.e.,  a compact hypersurface $\Sigma$ 
satisfies The Cauchy-Crofton Principle provided each point of $\Sigma$ 
has a neighborhood that does.  It then follows from Section 2 that:
 
  \proclaim {Proposition 2} To prove The Cauchy-Crofton Principle for a general compact, smooth hypersurface 
$\Sigma$ in $\Reals^n$ it suffices to treat the case of a smooth graph hypersurface.
\endproclaim

\proclaim {Continuity Lemma} If $\Sigma$ is a compact, smooth hypersurface of $\Reals^n$ then 
$f \mapsto \int_{\L^n} f^\Sigma(\ell) \, d\mu(l)$ is a well-defined continuous linear functional on the 
Banach space $L^\infty(\Sigma)$ of bounded, integrable, real-valued functions on $\Sigma$, i.e., 
it is continuous w.r.t.  the $\sup$ norm, $||\,||_\infty$.
\endproclaim

\proof
Recall that for $\ell \in \L^n$,  $\#(\ell \cap\Sigma)$ is defined to be zero if $\ell$ does not meet $\Sigma$ 
or if it meets it non-transversally, and otherwise it is the number of points in $\ell \cap\Sigma$. 
By Classic Cauchy-Crofton  $A(\Sigma) = {1\over \kappa_{n-1}} \int_{\L^n}  \#(\ell \cap\Sigma)\, d\mu(\ell)$.
For non-negative integers $k$ define 
\hbox{$\L^\Sigma_k := \{ \ell \in \L^n \mid \#(\ell \cap\Sigma) = k\}$}, so 
$\L^n$ is the disjoint union of the $\L^\Sigma_k$ and the set of measure zero where $\ell$ 
does not meet $\Sigma$ transversally, and hence  
$\sum_{k=1}^\infty k \mu(\L^\Sigma_k ) = \kappa_{n-1} A(\Sigma) $. 
If $f \in L^k(\Sigma)$, then by definition of $f^\Sigma$,
$|f^\Sigma(\ell)| \le k ||f||_\infty$ for $\ell \in \L^\Sigma_k $.
Hence :
$$\eqalign{
\left | \int_{\L^n} f^\Sigma(\ell) \, d\mu(\ell)\right | 
&= \left |\sum_{k=1}^\infty \int_{\L_k^\Sigma} f^\Sigma (\ell) \, d\mu(\ell) \right |\cr
&\le  \sum_{k=1}^\infty  \left |  \int_{L^\Sigma_k} f^\Sigma (\ell) \, d\mu(\ell) \right | \cr
&\le  \sum_{k=1}^\infty   \int_{L^\Sigma_k} \left | f^\Sigma (\ell) \right |  \, d\mu(\ell)  
\le  \sum_{k=1}^\infty k||f||_\infty \mu({L^\Sigma_k}) =  \kappa_{n-1} A(\Sigma) ||f||_\infty.     \cr
 }   $$ 
showing that $f \mapsto \int_{\L^n} f^\Sigma(\ell) \, d\mu(l)$  is bounded w.r.t. $||\,||_\infty$. \qed 

 Completing the proof of The Cauchy-Crofton Principle is now easy. By the usual argument it is enough 
 to prove it for a surface $\Sigma$ that is the graph of a continuous function $g: [0,1]^{n-1} \to \Reals$.
If a function $f : \Sigma \to \Reals$ is the characteristic function of the image under $g$ of a sub-cube of
$[0,1]^{n-1}$ then The Cauchy-Crofton Principle for $f$ follows from Classic Cauchy-Crofton, and so by 
the Linearity Lemma, The Cauchy-Crofton Principle also holds for ``step functions'', i.e., for linear 
combinations of such characteristic functions. But since any continuous $f$ is clearly a uniform limit 
of step functions, it follows from the Continuity Lemma that The Cauchy-Crofton Principle holds for 
all continuous $f: \Sigma \to \Reals$.  \qed

 \section{Constructing a Sequence $\{x_j\}$ Highly Equidistributed in a Hypersurface $\Sigma$}
 
   In this section  $\{\ell_j\}=\{(v_j,p_j)\}$  denotes a fixed \hbox{$k$-equidistributed} sequence in $\L_r^n, \ $ 
\footnote{\ for example, the one constructed in Subsection 4 of  Section 6.}  
and $\Sigma $ a bounded hypersurface included in $\B_r^n$.   We will show how to use 
the sequence  of lines $\ell_j$ to construct a \hbox{$k$-equidistributed} sequence $\{x_i$\} in 
$\Sigma$---essentially by just intersecting these lines with $\Sigma$. 

    We define a sequence $\sigma_j$ of non-negative integers by setting $\sigma_j$ to zero if $\ell_j$ 
does not meet $\Sigma$ or meets it non-transversely, and otherwise $\sigma_j$  is the positive integer 
denoting the cardinality of $\ell_j \cap \Sigma$; and we define $N_k := \sigma_1 + \ldots + \sigma_k$.  
(Thus, $N_k$ is the total number of points in which the first $k$ lines $\ell_j$ that are transversal 
to $\Sigma$ intersect $\Sigma$)
    
    To construct the sequence $\{x_i\}$, we start with the empty sequence 
 and consider each line $\ell_j$ in order and use it to append $\sigma_j$ 
 further points to the sequence $\{x_i\}$ as follows:
   \smallskip
{\narrower\noindent
 i) If $\ell_j$ is disjoint from $\Sigma$ or if it meets $\Sigma$ non-transversely (so that $\sigma_j = 0$)
then we append no new points to the sequence $\{x_i\}$. \par}
\smallskip
{\narrower\noindent
  ii) However, if $\ell_j$ does meet $\Sigma$, and meets it transversely, then it intersects
$\Sigma$ in the $\sigma_j$ points,  $\xi_i := \sigma^{(v_j,p_j)}(t_i)$,  
where $t_1 < t_2  <\ldots  < t_{\sigma_j}$, and in this case we append the $\sigma_j$ new points
$\xi_1, \ldots , \xi_{\sigma_j}$ of $\Sigma$ to the sequence $\{x_i\}$ in that order.       \par}

\smallskip\noindent
  Clearly, after processing the first $m$ lines in this way the sequence $\{x_i\}$ will have length $N_m$, 
and if $f: \Sigma \to \Reals$ is a bounded continuous function it is immediate from i) and ii) 
and the definition of  $f^\Sigma$ in section 8,  that 
$ \sum_{i=1}^m f^\Sigma(\ell_i) =  \sum_{j =1}^{N_{m} }f(x_j) $, so that:
$${1\over N_m} \sum_{j =1}^{N_{m} }f(x_j)  =  {\left(m \over N_m\right)} {1\over m}  \sum_{i=1}^m f^\Sigma(\ell_i) .$$

\proclaim{Lemma}  $$ \lim_{m \to \infty }  \left( {m \over\N_m} \right)  = { 1\over A(\Sigma) \,\kappa_{n-1} }.$$ 

\endproclaim

\proof
If we take $f \equiv 1$ in the Cauchy-Crofton Principle then, since $f^\Sigma(\ell)$ 
vanishes outside  $\L_r^n$, 
$A(\Sigma) = {1\over \kappa_{n-1}} \int_{\L_r^n} f^\Sigma(\ell)\, d\mu(\ell) $.
Since the sequence $\{\ell_j\}$ is by assumption $k$-equidistributed (and a fortiori
$1$-equidistributed) in $\L_r^n$,
$A(\Sigma) = {1\over \kappa_{n-1}}  \lim_{m \to \infty }{1\over m} (f^\Sigma(\ell_1) + \cdots + f^\Sigma(\ell_m))$.
But clearly  $f^\Sigma(\ell_j) = \sigma_j$ so $A(\Sigma) =  {1\over \kappa_{n-1}}\lim_{m \to \infty }{\N_m \over m}$.
\qed

  Applying the lemma and using again that the sequence $\{\ell_j\}$ is $1$-equidistributed it now follows that  
$$\lim_{m \to \infty} {1\over N_m} \sum_{j =1}^{N_{m} }f(x_j)  =  
\lim_{m \to \infty}  {\left(m \over N_m\right)} \lim_{m \to \infty} {1\over m}  \sum_{i=1}^m f^\Sigma(\ell_i)=
{ 1\over \kappa_{n-1} A(\Sigma)} \int_{\L^n} f^\Sigma(l)\, \,d\mu(l),$$
so by the Cauchy-Crofton Principle again, applied now to $f$, 
$$\lim_{m \to \infty} {1\over N_m} \sum_{j =1}^{N_{m} }f(x_j)  = {1\over A(\Sigma)} \int_\Sigma f(x)\, dx.$$

This shows that: 

\proclaim {Theorem} The sequence $\{x_j\}$ defined above satisfies 2) of Definition 1 of Section 3, and so is
$1$-equidistributed in $\Sigma$. 
\endproclaim  

To show that it  is $k$-equidistributed, repeat the above argument
but now for a bounded continuous function $f : \Sigma^k \to \Reals$ and instead of the Cauchy-Crofton Principle 
apply the  Multivariable Cauchy-Crofton Principle.  We omit these details.

\section{Principal Surface Types}

  Although we could continue developing the theory of point clouds living on hypersurfaces in 
$\Reals^n$, we will now specialize to the case $n=3$.  Not only does this make the discussion 
more intuitive and visualizable, but it also makes available some simplifications.

  We will start by describing the three principal types of surfaces that arise in 
mathematics and computer graphics,  {\it implicit surfaces\/},  {\it parametric surfaces\/}, 
and {\it triangulated surfaces\/}, discuss some of the important differences 
between them, and give a short intuitive description of the algorithms we use for generating
point clouds for each type. \smallskip

\item {a)} Implicit surfaces.
 
 An implicitly defined surface is defined as the level set of a function $f:\Reals^3 \to \Reals$.  
In more detail, a smooth implicit surface, $\Sigma$, is the set of solutions of an equation of 
the form $f(x_1,x_2,x_3)= 0$, where $f$ is a smooth real-valued function of three 
real variables. For example, the theorem of Pythagoras says that the sphere of radius 
$r$ centered at the origin consists of the points $(x_1,x_2,x_3)$ satisfying 
$x_1^2 + x_2^2 + x_3^2 = r^2$, and so it is an implicit surface 
with the defining function $f(x_1,x_2,x_3) = x_1^2 + x_2^2 + x_3^2 - r^2$.
For such surfaces it is easy to test if a point $(x_1,x_2,x_3)$ lies in the surface, since 
one only has to {\it evaluate\/} the function $f$ at the point---which we assume is easy---but it is 
often difficult to find points $(x_1,x_2,x_3)$ on the surface, since for that one has to {\it solve\/} 
the equation, and that can be a ``hard'' problem, even when $f$ is a polynomial.

We have already discussed in section 11 our algorithm, based on the  Cauchy-Crofton Formula, 
for generating a point cloud on an implicit surface. Intuitively speaking, the idea is to find the points 
where a set of lines (with highly equidistributed directions and closest points to the origin) intersect 
the implicit surface. We accomplish this by evaluating the function defining the surface along those 
lines and use a standard numerical technique to see where it vanishes (see Section 13 for details).

  A problem that must be addressed is that while some implicit surfaces, such as the 
sphere above, are bounded and have finite area, others (such as the plane, defined by 
$x_3 - x_1 - x_2 = 0$ or the paraboloid defined by $x_3 -x_1^2 -x_2^2 = 0$) are unbounded 
and have infinite area. For such surfaces, if the points in a point cloud are uniformly distributed 
in proportion to area, then the probability of a point being in a region of finite area 
would be zero! To avoid this problem we will, when constructing a point cloud,  
consider only a bounded piece of an implicit surface $\Sigma$; namely we replace
$\Sigma$ by $\Sigma_r$, its intersection with the disk $D^3_r$ of radius $r$ in $\Reals^3$, 
and construct point clouds lying entirely in $\Sigma_r$.
\medskip

\item {b)} Parametric surfaces. 
 
 A smooth parametric surface $\Sigma$ is defined by giving a smooth function:
 $$\Phi(u,v) := (X_1(u,v), X_2(u,v), X_3(u,v))$$
  from the $(u,v)$-plane into $(x_1,x_2,x_3)$-space.
 A point $(x_1,x_2,x_3)$ is  on $\Sigma$ if and only if there are values of the parameters 
$u$ and $v$ such that  $x_i = X_i(u,v)$ for $i = 1,2,3$.  The three real-valued functions $X_i$ 
that define the mapping $\Phi$ from $(u,v)$-space to $(x_1,x_2,x_3)$-space  are called the 
 parametric functions for the parametric surface. Usually we restrict the parameters 
 $(u,v)$  to lie in a restricted region of their plane (called the parameter domain) 
and often this is a rectangle with sides parallel to the axes.  For example,  the sphere 
of radius $r$ centered at the origin, considered above as an implicit surface, 
can also be considered as a parametric surface, with $u$ and $v$ 
respectively the latitude and longitude (measured in radians). The parameter 
domain is $-\pi/4 \le u \le \pi/4$ and  $-\pi/2 \le v \le \pi/2$ and $(X_1(u,v), X_2(u,v), X_3(u,v))$
are just the three spatial coordinates of the point of the sphere having latitude $u$ 
and longitude $v$.  (See {\it Remark 2\/} below for the explicit formulas.)
For parametric surfaces it is easy to find many points on the 
surface, just  substitute many pairs $(u,v)$ into $X_1,X_2,X_3$, but it is usually hard 
to decide if a given point $(x_1,x_2,x_3)$ lies on the surface, since this requires solving 
the three simultaneous equations $X_i(u,v) = x_i$  for $u$ and $v$. Our method for generating 
a point cloud for a parametric surface will be to define an associated triangulation of 
the surface and then use the method sketched below for triangulated surfaces.
\smallskip
\goodbreak

\item {c)} Triangulated surfaces. 
\smallskip

   Three distinct points, $v_1, v_2, v_3$,  in $\Reals^3$ define a planar surface 
called the {\it triangle\/} with {\it vertices\/}  $v_1, v_2$,  and $v_3$,  
and {\it edges\/} the line segments $[v_1,v_2]$, $[v_2,v_3]$, and $[v_3,v_1]$
joining pairs of vertices. The points of this triangle are just all points that lie on 
some line joining a vertex to a point on the opposite edge. Points 
of the triangle not on any edge are called {\it  interior\/} points of the triangle.
A triangulated surface $\Sigma$ is defined by giving a finite list of triangles
(called the {\it triangulation\/} of the surface) and $\Sigma$ is then 
just the union of these triangles.  Each vertex or edge of any triangle 
of the triangulation is called a vertex or edge of $\Sigma$, and 
if an edge of $\Sigma$ belongs to exactly one triangle of the triangulation it is called
a {\it boundary\/} edge, and its points are called boundary points of $\Sigma$.

Our method for generating a point cloud for a triangulated surface is based on Method 3 of 
Section 5.  We first compute the areas of the triangles to define a cumulative area distribution 
function. Each triangle is thereby associated with a subinterval of the unit interval whose length
is proportional to the area it contributes to the total  area of the surface. Then, given a random 
point in the unit interval, we can easily find which of these subintervals it belongs to, and then
place a random point in the associated triangle using the following algorithm.

 The standard $2$-simplex is given by $\Delta_2 := \{(u,v) \in [0,1]^2 \mid u + v \le 1\}$, and
we can define a pseudo-random sequence $\Rnd_{\Delta_2}(n)$ from $\Rnd_{[0,1)}$ 
using rejection; i.e., we define a sequence $\{ (u_k,v_k)\}$ in $[0,1]^2$ by 
$u_k :=  \Rnd_{[0,1)}(2k -1)$ and 
$v_k  := \Rnd_{[0,1)}(2k)$ and let  $\Rnd_{\Delta_2}(n)$ to be the $n$-th element of  
this sequence such that $u_n + v_n \le 1$. 

  Then if $\Delta$ is a triangle in $\Reals^3$ having vertices $v_1,v_2,v_3$, 
we can define a pseudo-random sequence for $\Delta$ by
$\Rnd_\Delta(n) := \Psi \circ \Rnd_{\Delta_2}(n)$, where $\Psi : \Delta_2 \to \Reals^3$
is the barycentric parametric representation $(u,v) \mapsto u v_1 + v v_2 + (1-u-v) v_3$
of  $\Delta$ (see Remark 3 above). This works 
because $\Psi$ is clearly an affine map and so multiplies all areas by a fixed
scale factor (its constant Jacobian determinant).

\medskip 
\goodbreak

\remark{Remark 1}
 For each of the above surfaces types, one needs extra  conditions to insure 
that what results is really what mathematicians call a surface (with boundary). 
\medskip
\item {a)}For an implicit surface, $\Sigma$, given by $f(x_1,x_2,x_3) = 0$, one needs to assume 
that the gradient of $f$ does not vanish anywhere on the surface, i.e., at any point 
$(x_1,x_2,x_3)$ of $\Sigma$, at least one of the three partial derivatives of $f$ 
is non-zero. This insures that $\Sigma$ has a two-dimensional tangent space at the point.  
In fact it then follows (from the Implicit Function Theorem) that near each of its points $\Sigma$ 
can be represented locally as the graph of a smooth function.
\medskip
\item {b)} For a parametric surface, $\Sigma$, given by  $\Phi(u,v)$, one needs to assume that if
$u  = \xi(t),\ v =\eta(t))$ are the equations of a straight line in the parameter domain, then
the image curve on the surface (given by $\Phi(\xi(t),\eta(t))$) has a non-zero tangent vector 
for all $t$. This again  insures that $\Sigma$ has a two-dimensional tangent space at every point.
\medskip
\item {c)}To make the definition of a triangulated surface $\Sigma$ rigorous, one needs the 
following three conditions:

\itemitem {i)} A point of $\Sigma$ that is an interior point of some triangle 
                    of the triangulation does not belong to any other triangle of the triangulation. 
\itemitem {ii)} Any non-boundary edge of $\Sigma$ should be an edge of exactly two 
                     triangles of the triangulation.
\itemitem {iii)} If a vertex $v$ of $\Sigma$ lies on an edge $e$, then $v$ is one of 
the two endpoints of $e$, and if $e$ is a boundary edge then $v$ is on exactly one 
other boundary edge.

\noindent
(The first two condition insure that $\Sigma$ looks right 
(i.e., in math jargon, ``has the right local topology'') at interior points and edges 
and the third condition insures that it looks right at all vertices.) However, the conditions 
a), b), c) do {\it not\/} enter into the algorithms for defining point clouds and so 
will not play a significant role in what follows.

\endremark
\medskip

\remark{Remark 2}
Sometimes an implicit surface can also be given parametrically, by explicitly solving 
its defining equation $f(x_1,x_2,x_3)= 0$.  For example a plane is given implicitly by an equation 
of the form $a_1x_1 + a_2 x_2 + a_3x_3 + a_4= 0$, or (assuming $a_3 \ne 0$)  parametrically by 
$x_1 = u,\  x_2 = v, \ x_3 = ( a_4 - a_1 u -a_2 v)/a_3$, 
and similarly the implicitly defined ellipsoid 
$(x_1/a_1)^2 + (x_2/a_2)^2 +(x_3/a_3)^2 - 1 =0$ 
can be given parametrically by 
$x_1 = a_1 \cos(u) \sin(v), \ x_2 = a_2 \sin(u)\sin(v), x_3 = a_3\cos(v)$. 
But these are exceptional cases, and in general no such explicit solution 
of an implicit equation in terms of elementary functions is possible.
\endremark
\medskip

\remark{Remark 3}
A triangle $\Delta$ with vertices $v_1, v_2$, $v_3$ can be considered 
a parametric surface with parameter domain $\Delta_2$, the {\it standard $2$-simplex\/},
consisting of  all $(u,v) \in [0,1]^2$ with $u + v \le 1$, and defining 
$\Psi(u,v) := u v_1 + v v_2 + (1-u-v) v_3$. We call $\Psi(u,v)$ the point of 
$\Delta$ with {\it barycentric coordinates\/}  $u,v$, and $1-u-v$. 
Thus a triangulated surface can be considered as a finite union of parametric surfaces. 
In the reverse direction, it is easy to {\it approximate\/} a parametric surface $\Sigma$ by a 
triangulated surface $\Sigma'$. If we first triangulate the parameter domain of $\Sigma$ into 
small triangles, and then take the images of the vertices of these triangles as 
the vertices of triangles of a triangulation, we get a triangulated surface $\Sigma'$ 
that approximates $\Sigma$. 
\endremark
\medskip
\goodbreak

\remark{Remark 4}
While the surfaces that occur in mathematical contexts  are usually of the 
 implicit or parametric type, the surfaces of interest in computer graphics 
 are frequently boundaries of some real world solid object, 
 such as a human face or a teapot. Such surfaces are in general 
 too irregular to have convenient representations as parametric 
 or implicit surfaces, so it has been customary to represent them 
 by ad hoc triangulations.  However, recently there has been rapid progress 
 in a process called LIDAR or laser scanning, a method for the rapid construction 
 of a  ``point cloud''
\footnote{The points of these laser generated point clouds are {\bf not} uniformly 
 distributed or independent.}
 from a solid object. It works by scanning a rapidly pulsed laser beam over the 
 surface of the object and collecting the reflected pulses in several optical devices. 
 The 3D location of a point on the surface from which a pulse is reflected 
 is calculated either by triangulation or by time-of-flight methods, and the
 totality of these calculated 3D locations gives the point cloud.
 There is associated software that can take the resulting LIDAR point cloud 
 and construct from it a so-called Voronoi decomposition of the surface, and 
 then from that in turn construct a triangulation. So it could be 
 argued that the point clouds created by laser scanning are merely 
 steppingstones to triangulations. But in fact, a major reason for creating 
 the triangulation is that so much software has been developed 
 over the years to work with triangulated surfaces, while point clouds
 are relatively new. Since there is a loss of information in going from 
 point cloud to triangulation, it seems likely that the point cloud data 
 structure will gradually become another primary means for storing and 
 working with surfaces---and perhaps eventually even the predominant one. 
This will make it important to also have good algorithms for creating high 
quality, uniform point clouds from implicit, parametric, and triangulated surfaces, 
in order to be able to work with all the various types of surfaces in a uniform 
way. Indeed, in virtual reality applications, one already sees combinations of 
real-world objects and mathematically defined surfaces in the same scene.
It will also be important to be able to take a non-uniform LIDAR generated 
point cloud and ``upgrade'' it to one that is uniform. One way to do that is to 
construct a triangulated surface from the point cloud as outlined above and 
apply our algorithm for creating a uniform point cloud from a triangulation.
\endremark
\smallskip

\remark{Remark 5}
In addition to the somewhat conjectural reason just alluded to, there are many  
current applications requiring a point clouds approach. One is Monte Carlo 
based numerical methods, and we put off discussion of this until later. 
Another is visualization: when mathematicians would like to examine the geometry 
of a surface, either for research or teaching purposes,  they can now use mathematical 
visualization software to create and rotate an image of the surface on a computer screen.  
As discussed earlier, the classic tool for displaying implicit surfaces is 
raytracing, and while this method produces realistic renderings, 
it does have drawbacks. In a raytrace of a surface with 
many layers, the layer nearest the viewpoint obscures those further away, and 
since raytracing is a comparatively slow process, it is difficult to rotate 
a raytraced image in real time to see the surface from different perspectives.  
Point cloud rendering is free of these drawbacks. 
  \endremark

\section{Point Clouds for Implicit Surfaces}

\subsection{Intersecting a line with an implicit surface}

  As remarked earlier,  it may be difficult to find solutions of  
$f(x_1,x_2,x_3) = 0$ in a rectangular region $m_i \le x_i \le M_i$. 
However, for the analogous single variable problem,
namely that of solving the equation $g(t) = 0$ 
on an interval $m \le t \le M$,  there  {\bf is} an effective 
numerical approach to finding solutions. 
Recall that the Intermediate Value Theorem tells us that 
a continuous real-valued function of a real variable that has    
opposite signs at the endpoints of some interval vanishes  
for some point in that interval. So a good numerical method for finding zeros of 
$g(t)$ on $m \le x \le M$ is to divide $[m,M]$ into small  
subintervals, check on which subintervals $g$ changes sign, and then use 
either bisection or Regula Falsi to locate quickly a zero on each of 
these subintervals.
 If $x(t) =  t v +p$ is the equation of a straight line $L$ in $\Reals^3$ 
and we define $g(t) := f( x(t))$, we can use this method  
to find the solutions of $g(t) = 0$, i.e., to find the intersections of the line 
$L$ with the implicit surface $\Sigma$ defined by the equation $f = 0$, and 
this will be our point generating method for creating a point cloud on $\Sigma$ 
(or rather on $\Sigma_r$, the part of $\Sigma$ that lies inside the ball of radius $r$).

But how should we choose the lines $L$ to insure that the points obtained 
in this way are uniformly distributed? We consider this next.
\medskip\goodbreak

\subsection{Point clouds for Implicit Surfaces} 

We need an algorithm $\Rnd_{\L^3_r}$ for selecting lines $(v,p)$
in $\L^3_r$ uniformly with respect to kinematic measure.  It is clear 
from the definition of kinematic measure that we get such a selection 
algorithm by first selecting $v \in \S^2$ using $\Rnd_{\S^2}$, 
and then selecting $p \in v^\perp$ using the following algorithm 
$\Rnd_{D^2_r(v^\perp)}$ for selecting a pseudo-random sequence of 
points in $D^2_r(v^\perp)$,  the disk of radius $r$ in $v^\perp$.

Let $J_{2,3}(x_1,x_2) := (x_1,x_2,0)$ be the inclusion of $\Reals^2$ into $\Reals^3$
and $T_v$ the rotation of $\Reals^3$ that is the identity on the line $P^\perp$ 
orthogonal to the plane $P$ spanned by $v$ and $e_3 := (0,0,1)$ and on 
$P$  rotates $e_3$ to $v$.   Since $T_v$ maps $e_3$ to $v$, it maps 
$e_3^\perp$, the image of $J_{23}$, onto $v^\perp$. Then we define 
$\Rnd_{D^2_S(v^\perp)}(n) := T_v \circ J_{2,3} \circ \Rnd_{D^2_S}(n)$,
i.e., we use $J_{23}$ to inject $\Rnd_{D^2_S}$   into $e_3^\perp$, 
and then rotate it with $T_v$ into $v^\perp$. 
 $T_v$ is a transvection, i.e.,  if $b$ is the midpoint  
of the great circle joining $e_3$ to $v$, then $T_v$ is the composition of reflections 
in two one-dimensional subspaces, namely reflection $\rho_b$ in 
the line spanned by $b$ followed by reflection $\rho_v$ in the line 
spanned by $v$, so finally our algorithm reads:  
$\Rnd_{D^2_S(v^\perp)}(n) := \rho_v\circ\rho_b \circ J_{2,3} \circ \Rnd_{D^2_S}(n)$.

It now follows from the corollary of the Generalized Cauchy-Crofton Formula 
that if $\Sigma$ is an implicit surface, then the expected number of 
intersections such a pseudo-randomly selected line will have with 
any region of $\Sigma_r$ is proportional to the area of that region. 
Thus, selecting a sequence of lines with $\Rnd_{\L^3_r}$ and then finding the intersections
of these line with $\Sigma_r$, using the method just described above, 
provides our sought for algorithm for constructing a uniformly distributed sequence
for an implicit surface $\Sigma$.

\section{Point Clouds for Triangulated Surfaces}

  Let $\Sigma$ be a triangulated surface with triangulation given by a list of triangles 
$\Delta_i$, $i = 1,2,\ldots, N$, having vertices $v^i_1,v^i_2,v^i_3$. We recall that 
 the area of $\Delta_i$ is given by
 $A(\Delta_i) := {1\over 2}||(v^i_2 -v^i_1) \times  |(v^i_3 -v^i_2) ||$.
We define $K_0 := 0$ and $K_i  := K_{i-1} + A(\Delta_i)$, 
$i = 1,\ldots, N$;  so $K_i$ is the sum of the areas of triangles $\Delta_1$ 
through $\Delta_i$ and in particular $K_N = A(\Sigma)$, the total area of $\Sigma$.
Thus  if we define $\delta_i := {K_i\over K_N}$, then the sequence 
$0=\delta_0  < \delta_1 < \ldots < \delta_N = 1$ defines a partition of $[0,1]$ 
into $N$ subintervals  $I_j := [\delta_{j-1}, \delta_j]$ with the length of the $j$-th 
subinterval being the fraction of the total area of $\Sigma$ that is in $\Delta_j$.
This provides just what we need to define an algorithm, $\Rnd_\Sigma$, for choosing 
uniformly distributed points of $\Sigma$. 
Namely, select a pseudo-random number $x =\Rnd_{[0,1)}$ in $[0,1]$, 
determine the subinterval $I_j$ that  this $x$ belongs to, and let $\Rnd_\Sigma := \Rnd_{\Delta_j}$, 
where  $\Rnd_{\Delta_j}$ denotes the pseudo-random selection function for $\Delta_j$ defined 
in  12 c).

  A straightforward way to find the index $j$ such that $x$ belongs to $I_j$ is to 
test the inequality $\delta_i \le x$ for $i = 1,2,\ldots$ and let $j$ be the first index $i$
for which it is true. However this takes on the order of $N$ steps, and since in 
practice $N$ will be on the order of $10^4$ to $10^5$, it is much better to use a 
bisection search, which will produce $j$ in an order of $\log_2 N$ steps.

\section{Point Clouds for Parametric Surfaces}

  Let $\Phi: [uMin,uMax] \times [vMin,vMax] \to \Reals^3$ define a parametric surface $\Sigma$.
We will refer to its domain as the {\it parameter rectangle\/}.  A standard 
technique in computer graphics is to divide the parameter rectangle into a grid of
sub-rectangles by dividing each of the two intervals $[uMin,uMax]$ and 
$[vMin,vMax]$ intervals into equal length sub-intervals. By dividing each 
sub-rectangle by a diagonal into two triangles and mapping these triangles 
with $\Phi$, we get a triangulation of a piecewise linear approximation $\Sigma'$ 
of $\Sigma$. Thus the point cloud algorithm for $\Sigma'$ in the preceding section will
give an approximate point cloud for $\Sigma$, and by a slight modification (see below) 
we can also get a point cloud for $\Sigma$ itself.

\smallskip\noindent
In more detail, let $uRes >1$ and $vRes >1$ be integers. We divide 
$[uMin,uMax]$ into the $(uRes - 1)$ subintervals $[u_i, u_{i+1}]$, $i = 1,\ldots uRes -1$, 
where we define $uStep := (uMax - uMin)/(uRes - 1)$ and 
$u_i := uMin  + (i-1)*uStep$, and in a similar manner we divide $[vMin,vMax]$ 
into the $(vRes - 1)$ subintervals $[v_i, v_{i+1}]$, $i = 1,\ldots vRes -1$. 
The point $g_{i,j} := (u_i,v_j)$ of the parameter rectangle are called the parameter grid-points
and their images $\gamma_{i,j} := \Phi(g_{i,j})$ are called the surface grid points.
The rectangles $r_{i,j}$ with vertices $(g_{i,j},g_{i+1,j},g_{i+1,j+1},g_{i,j+1})$ are called 
the parameter sub-rectangles and we get a triangulation of the parameter 
rectangle by dividing each $r_{i,j}$ into the two triangles 
$\delta^+_{i,j} := (g_{i,j},g_{i+1,j+1},g_{i,j+1})$ and 
$\delta^-_{i,j} := (g_{i,j},g_{i+1,j},g_{i+1,j+1})$.  Finally we  get 
a triangulated surface $\Sigma'$ by taking their images 
$\Delta^+_{i,j} := \Phi(\delta^+_{i,j})= (\gamma_{i,j},\gamma_{i+1,j+1},\gamma_{i,j+1})$ 
and
$\Delta^-_{i,j} := \Phi(\delta^-_{i,j})= (\gamma_{i,j},\gamma_{i+1,j},\gamma_{i+1,j+1})$. 

 Note that the surface grid points $\gamma_{i,j}$ belong to the surface $\Sigma$ by
their definition, so the vertices of the triangles $\Delta^+_{i,j}$ and $\Delta^-_{i,j}$ lie 
in $\Sigma$. Thus the smooth surfaces $\Sigma$ and the piecewise linear surface 
$\Sigma'$ coincide at the surface grid points, and it follows from the differentiability 
of $\Phi$ that, as $uRes$ and $vRes$ approach infinity, the surface $\Sigma'$ 
converges uniformly to $\Sigma$. Hence, if we construct a point cloud for the 
triangulated surface $\Sigma'$, using the algorithm of the preceding section, then 
for large $uRes$ and $vRes$ this will be a good approximation of a point cloud 
for $\Sigma$, although the points will only lie near $\Sigma$ not on it.  However it
is easy to modify the algorithm for the triangulation point cloud slightly so that the 
points do lie on $\Sigma$. Namely in the definition of the selection function 
$\Rnd_\Delta$, replace the barycentric parameterization $\Psi$ by $\Phi$.

\section{Point Cloud Surface Normals}

   In many applications of point clouds associated to a surface $\Sigma$, it is important 
to know the unit vector $\nu_p$ normal to $\Sigma$ at points $p$ of the cloud.
One important application of $\nu_p$ arises when we try to visualize the surface by
rendering the point cloud on a computer screen. After projecting each point to the 
view plane representing the screen we could simply render these points black on a
white background. However, for a more realistic looking version of the surface, 
we  can``paint'' each point of the cloud the color it would get in a raytrace of the 
surface or,  even better, use that color to paint a small disk tangent to the surface at the 
point. In fact, if the point cloud is sufficiently dense, such a rendering provides a quite good 
approximate raytrace that is much faster than a full one.  Now to calculate this raytrace 
color at a cloud point $p$ we have to know how a light ray from each of several colored 
light sources reflects off the surface when striking it at $p$, and by the Law of Reflection 
for specular surfaces (``Angle of Incidence Equals Angle of Reflection''), this depends in 
a crucial way on $\nu_p$. And of course the tangent plane at $p$, $T\Sigma_p= \nu_p^\perp$, 
is also determined by $\nu_p$.

  Other uses of $\nu_p$ follow from the fact that it can be used  to 
approximate $\Sigma$ near $p$. Indeed,  by choosing two vectors $f_1$ and $f_2$ 
that are linearly independent of $\nu_p$ and applying the Gram-Schmidt algorithm to
 $\nu_p, f_1, f_2$, we get an orthonormal basis $\e_1, \e_2$ for $T\Sigma_p$, so that 
$\Phi(u,v) := u \e_1 + v \e_2$  gives a parameterization of $T\Sigma_p$ that 
by Taylor's Theorem will be a good first order approximation to $\Sigma$ near $p$.

   For the point cloud algorithms for the three types of surfaces $\Sigma$ discussed 
above, it is easy to construct the normal at a cloud point $p$ using standard methods 
of multi-variable Calculus. For example, if $\Sigma$ is given by the implicit 
equation $f(x_1,x_2,x_3) = 0$, we can first calculate $\nabla f_p$ by computing the 
three partial derivatives of $f$, and then normalize it to get $\nu_p$.

  But what if we are given {\it only\/} the point cloud itself and do not know the surface $\Sigma$ 
that it came from---for example, suppose the point cloud comes from a LIDAR scan of a solid 
object. Is there some way we can nevertheless construct at least a good approximation to
the normal $\nu_p$ to $\Sigma$ at a cloud point $p$? If the point cloud is dense enough 
(relative to the maximal curvature of $\Sigma$) then there is a fairly obvious approach:
choose two other points $q$ and $r$ of the cloud that are close to $p$ and normalize the 
cross-product $(q-p) \times (r-p)$. And, to get a more accurate approximation to the normal, 
one can average the cross-products $(q_i-p) \times (r_i-p)$ for a number of different 
pairs of close neighbors $(q_i,r_i)$ of $p$.  It might look like there is a catch in this; 
finding those good close neighbors of $p$ seems to require considerable computation, 
involving sorting the points of the cloud in various ways. But fortunately it turns out 
to be possible to carry it out surprisingly fast. 
  
  \vskip 15 pt

\section{Applications: What are Point Clouds Good for?}

\subsection{Mathematical Visualization}

   In the mid-1990s, while developing our mathematical visualization program 
(3D-XplorMath), we became dissatisfied with raytrace based methods for visualizing 
implicit surfaces, and our interest in point clouds grew out of attempts to improve on the 
tools available for this purpose.  CPUs and  graphic chips have improved considerably in speed 
since that time, and today a basic raytrace takes only a fraction of a minute on a modest laptop 
computer, but back then that process seemed glacially slow, even on a high-end desktop machine. 
In addition, for an immersed implicit surface that has many overlapping layers, the layer closest 
to the viewing camera can obscure interesting features on parts of the surface  behind it, 
\footnote{This problem can be at least partially overcome by using ``transparency'', however
doing so slows down a raytrace considerably. }
and even today the time required to redraw each frame is too slow to make mouse 
rotation of the surface an acceptable cure for that problem.   \par
 
  Point cloud rendering of implicit surfaces provided our sought for solution to these 
 problems, particularly when coupled with the anaglyph  
(i.e., red/blue glasses based) stereo viewing technique. It is quite striking to see simultaneously  
the many layers of a complex surface floating in front of you, almost as if the surface were constructed 
out of vanishingly thin glass---and it becomes even more remarkable when this image is set smoothly rotating.
Indeed, point cloud rendering performed so well for seeing implicit surfaces that we ended up including 
point cloud based techniques for viewing parametric and triangulated surfaces as a supplement  
to the more classic rendering methods already available for those surfaces.
   \par

   In our original approach to generating implicit surface point clouds we
intersected the surface with randomly chosen lines, each of which
was parallel to one of the three coordinate axes. While such point clouds are
reasonably uniform---enough so for most visualization purposes, 
a computation shows that their density can vary by a factor $\sqrt3$
from one point to another on a surface, and this renders them unsuitable for use with 
Monte Carlo methods (see below) or other numerical applications whose accuracy 
depends on having highly uniform distributions of points. After considerable reflection 
and experimentation we realized that the Cauchy-Crofton Formula implied that if we 
selected the lines used to intersect the surface to be uniformly distributed with respect to 
kinematic measure, the resulting point clouds would have the uniformity properties required 
to support such computational applications.  There is also a subtle visualization benefit 
from point cloud uniformity;  namely as one looks towards a point on a surface ``contour'', 
i.e., a point where the line of sight is tangent to the surface, the points of the cloud appear 
increasingly dense, and this not only helps make the contours visible but also gives the surface 
a more accurate three dimensional appearance.
\par

  \subsection{Monte Carlo Methods} 

It has been said that Monte Carlo methods ``Lift the Curse of Dimensionality.''  Let's see 
what that means.    Suppose we have a $C^1$ function defined on the $n$-cube,
$f: I^n \to \Reals$,
and we want to estimate its integral $J =\int f(x) \,dx$. A natural measure 
of the complexity of a numerical algorithm for estimating $J$ is the number $N(\eps)$ 
of evaluations of $f$ required by the algorithm in order that the expected error in the  
estimate $J^\prime$ it generates will differ from the correct value by less than $\eps$. 
We will estimate this is for both a simple Riemann integral approach and using a Monte Carlo 
method.  For the Riemann integral inspired calculation we divide each side of $I^n$ into $k$ equal
sub-intervals of length $1\over k$, thereby partitioning the cube into $N_k :=k^n$ sub-cubes 
$\kappa_i$ of diameter $\delta_k = {\sqrt{n}\over k}$ and volume $v_k = {1\over N_k}$. 
To estimate $J$, we select arbitrary points $x_i$, one  in each $\kappa_i$, and compute the
Riemann  sum $J_k := \sum_{i=1}^{N_k} f(x_i) v_k $, and we note that this requires $N_k = k^n$ 
evaluations of $f$.  On the other hand, if 
$\bar x_i$ in $\kappa_i$ is a point where $f$ assumes its mean value in $\kappa_i$,
(i.e., such that $f(\bar x_i) = {1\over v_k} \int_{\kappa_i} f(x) \, dx$) then 
$J= \sum_{i=1}^{N_k} f(\bar x_i) v_k $, and hence 
$|J -J_k| \le \sum_{i=1}^{N_k}| f(\bar x_i) - f(x_i)| v_k $.
Since $f$ is $C^1$, $K := \sup ||{\nabla f}_x||$ is a Lipshcitz constant 
for $f$, hence if $||\bar x_i - x_i || < {\eps \over K}$ for all $i$ then 
$ | f(\bar x_i) - f(x_i)| < \eps$ for all $i$ and hence $|J -J_k| \le \eps$. 
Since $x_i$ and $\bar x_i$ are arbitrary elements of $\kappa_i$, to insure 
$||\bar x_i - x_i || < {\eps \over K}$ we need the diameter  $\delta_k$ of the $\kappa_i$ 
less than ${\eps \over K}$, i.e., ${\sqrt{n}\over k} < {\eps \over K}$ and hence 
$k > {K \sqrt{n} \over \eps}$. This means that for this algorithm the measure of 
complexity $N(\eps) = N_k = k^n \ge \left({{K \sqrt{n} \over \eps}}\right)^n$ and we se that, 
using standard ``Big O'' notation, for this algorithm $N(\eps) = {\bold O}(\eps ^{-n})$.  
While this estimate may at first glance appear innocuous, on further examination it  
leads us to typical examples of the ``Curse of Dimensionality''. 
For there are many examples in theoretical physics where the answer to an important 
question requires computing a definite integral over $I^n$ where $n > 100$. Suppose we are 
satisfied with an accuracy of only $0.01$ in our calculation. Then the above method for evaluating 
the integral  is clearly hopeless, since it requires on the order of $100^{100}$ evaluations of $f$! 
Another way to look at the curse is to note that it says if one needs $N$ evaluations to get the 
expected error less than $\eps$, then it requires $2^n N$ evaluations to get it less than   $\eps \over 2$.
While it is true that similar but more sophisticated numerical methods, such as Simpson's Rule,
can cut down somewhat on the order of growth of $N(\eps)$ as a function of the dimension $n$,
this improvement is too small to solve the problem---the curse remains.  What we need is a method 
for which the order $N(\eps)$ of the error estimate is independent of $n$, and the remarkable fact 
is that the Monte Carlo approach provides just that. We will not prove it here (for details, see [J] p. 291) 
but it follows from the Central Limit Theorem that if we use an $n$-distributed RNG to estimate 
the integral using the Monte Carlo formula of Definition 2 in Section 3, then the expected error with $N$ 
evaluations is ${\bold O}\left({1\over \sqrt{N}}\right)$ so $N(\eps) = \bold O(\eps ^{-2})$, and we only require 
four times as many function evaluations to double the expected accuracy.
  
\section{Generalizations}

  While the theory developed above is restricted to hypersurfaces, it generalizes to provide a method for 
constructing equidistributed point clouds on submanifolds $\Sigma$ of $\Reals^n$ of higher codimension, 
say $k$. Although the details get somewhat more complex,  the main constructions and proofs go over in 
a straightforward manner. The first modification is that the space $\L^n$ of affine lines in $\Reals^n$ 
gets replaced by the space $\Aff_k^n$ of $k$-dimensional affine subspaces of $\Reals^n$. 
(Of course, $\Aff_1^n = \L^n$.)  The group of Euclidean motions acts transitively on 
$\Aff_k^n$ and again there is a natural kinematic measure 
$d\mu$ that is invariant under this action. Except for a set of $V \in \Aff_k^n$ having kinematic 
measure zero, the number $\#(V \cap \Sigma)$ of points where $V$ intersects $\Sigma$ is finite, 
and the obvious generalization of Cauchy-Crofton is valid, namely 
$\int_{\Aff_k^n} \#(V \cap \Sigma) \, d\mu(V)$,
the average number of intersection points, is a universal constant times the volume of $\Sigma$. 
When it comes to actually implementing this generalization numerically, there appears 
to be a problem; the technique described in the first subsection of Section 13 for calculating where 
a line meets an implicit hypersurface was based on the intermediate value theorem, i.e., the principle 
that if a continuous real-valued function has different signs at the endpoints of an interval 
then it vanishes somewhere in the interval, and this requires $k = 1$. Now even when $k=1$, 
once we have a good approximation to the intersection point, it is more efficient to use Newton's Method 
rather than bisection to find it with precision,  and it is a higher dimensional version of this Newton's Method 
approach that saves the day. The fact that Newton's Method generalizes to an effective algorithm 
for solving our problem goes back to an old paper of Kantorovitch [Ka] and it is usually known as the 
Newton-Kantorovich Theorem. We will not go into further detail here and instead refer the interested 
reader to a nice exposition by J. M. Ortega [O].
What makes this generalization to higher codimensional submanifolds of significant interest is 
 that, coupled with the Moser-Whitney Theorem (section 1.5), it extends our methods for
 constructing highly equidistributed point clouds to arbitrary compact measured manifolds, 
 and thereby provides an approach to implementing Monte Carlo methods in that degree of 
 generality.  We have not tried to estimate the complexity of the combination of algorithms that 
 are required for such an implementation and consider that to be a worthwhile project.

\bigskip\bigskip
 

\centerline{\bf References}
\medskip

\item{[B]} {Billingsley, P., {\it Convergence of Probability Measures,\/}  2nd Ed., John Wiley, 1999 }
\medskip
\item{[Bl]} {Blaschke, W., {\it Vorlesungen \"uber Integralgeometrie,  Bd.I 1936; Bd.II 1937 } Teubner Verlag, Leipzig} 
\medskip
\item{[B-C]} {Bailey, C.H. and Crandal, R., {\it Random Generators and Normal Numbers}, Experiment. Math., A.K. Peters, 2002, pp. 527--546} 
\medskip
\item{[B-T]} {Bott, R. and Tu, L.., {\it Differential Forms in Algebraic Topology}, Springer Verlag, 1982} 
\medskip
\item{[C]} {Chern, S.S., {\it On Integral Geometry in Klein Spaces}, Ann. of Math., {\bf 43}(2), 1942, pp. 178--189} 
\medskip
\item{[Cr]} {Crofton, M.W., {\it Probability}, Encyclopedia Britannica, 9th Ed.,  vol. 19, 1885, pp.768--788} 
\medskip
\item{[E]} {L'Ecuyer, P., {\it Uniform Random Number Generators: A Review}, 
Proceedings of the I997 Winter Simulation Conference, 1997, pp.127--134} 
\medskip
\item{[F]} {Franklin, J.N., {\it Deterministic Simulation of Random Processes}, Math. Computation, {\bf 17}, 1963, pp. 28--57} 
\medskip
\item{[Fe]} {Federer, H., {\it Colloquium Lectures on Geometric Measure Theory}, Bull. Amer. Math. Soc., {\bf 84}, 1978, pp. 291--338} 
\medskip
\item{[F-W1]} {Fuselier, E. and Wright,G., {\it A High-Order Kernel Method for Diffusion and Reaction-Diffusion Equations on Surfaces}, 
Journal of Scientific Computing, {\bf 56}(3), Springer, 2013, pp. 535--565} 
\medskip
\item{[F-W2]} {Fuselier, E. and Wright,G., {\it Scattered Data Interpolation on Embedded Submanifolds with Restricted 
Positive Definite Kernels: Sobolev Error Estimates}, SIAM Journal on Numerical Analysis, {\bf 50}(3), 2012, pp. 1753--1776} 
\medskip
\item{[H]} {Hirsch, M.W., {\it Differential Topology}, Springer Verlag, 1976} 
\medskip
\item{[J]} { Judd, K.L., {\it Numerical Methods in Economics}, MIT Press, 1998} 
\medskip
\item{[Ka]} {Kantorovich, L. V., {\it On NewtonÕs method for functional equations}, Dokl. Akad. Nauk. SSSR, {\bf 59}, 1949, pp. 1237--1240} 
\medskip
\item{[K1]} {Knuth, D., {\it The Art of Computer Programming, Volume 2:   Seminumerical Algorithms,\/}  Addison-Wesley, 1997}
\medskip
\item{[K2]} {Knuth, D.,  {\it Construction of a Random Sequence}, BIT, {\bf 5}, 1965, pp. 246--250} 
\medskip
\item{[K-R]} {Klain, D. and Rota G.-C., {\it Introduction to Geometric Probability,\/}  Cambridge Univ. Press, Cambridge, 1997}
\medskip
\item{[K-N]} {Kuipers, L. and Niederreiter, H., {\it Uniform Distribution of Sequences,\/}  Wiley, New York, 1974}
\medskip
\item{[L-Z]} {Liang, J. and Zhao, H., {\it Solving partial differential equations on point clouds}, 
SIAM Journal on Scientific Computing, {\bf 35}(3), 2013, pp. 1461--1486} 
\medskip
\item{[M]} {Moser, J., {\it On the Volume Elements on a Manifold}, 
Transactions of the AMS, {\bf 120}(2), 1965, pp. 286--294} 
\medskip
\item{[M-N]} {Matsumoto, M.; Nishimura, T., {\it Mersenne twister: a 623-dimensionally equidistributed 
uniform pseudo-random number generator\/}, 
ACM Transactions on Modeling and Computer Simulation, {\bf 8}, 1998, pp. 3--30} 
\medskip
\item{[O]} {Ortega, J.M., {\it The Newton-Kantorovich Theorem\/}, The American Math. Monthly, {\bf 75}(6), 1968, pp. 658--660} 
\medskip
\item{[P]} {Palais, R., {\it The Visualization of Mathematics: Towards a Mathematical Exploratorium\/}, 
Notices of the AMS, {\bf 46}(6), 1999, pp. 647--658} 
\medskip
\item{[P-F]} {Paiva, J.C. \'Alvarez and  Fernandes, E., {\it Gelfand Transforms and Crofton Formulas\/}, 
Sel. Math., New Ser., {\bf 13}, 2007, pp. 369--390} 
\medskip
\item{[R]} {Riesz, F., {\it Sur les op\'erations fonctionnelles lin\'eaires\/}, 
C.R. Acad. Sci., Paris, {\bf 149}, 1909, pp.  974--977 } 
\medskip
\item{[S-1]} {Santalo, L., {\it Introduction to Integral Geometry,\/}  Hermann, Paris, 1953}
\medskip
\item{[S-2]} {Santalo, L., {\it Integral Geometry and Geometric Probability,\/}  Addison-Wesley, 1976}
\medskip
\item{[SP1]} {Spivak, M., {\it Calculus on Manifolds}, W.A. Benjamin, 1965} 
\medskip
\item{[SP2]} {Spivak, M., {\it A Comprehensive Introduction to Differential Geometry Vol. I}, Publish or Perish, 1999} 
\medskip
\item{[U]} {Metropolis, N. and Ulam, S., {\it The Monte Carlo method\/}, Journal of the American Statistical Association, {\bf 44}, 1949, pp. 335--341} 
\medskip
\item{[VN]} {von Neumann, J., {\it Various techniques used in connection with random digits\/}, 
Nat. Bur. Stand., Appl. Math. Series {\bf 12}(1), 1951, pp. 36-88} 
(Reprinted in von Neumann's Collected Works, 5 (1963), Pergamon Press pp 768-770)
\medskip
\item{[Wh]} {Whitney H., {\it Differentiable manifolds\/}. Ann. of Math. {\bf 37}(2)  (1936), no.3, pp. 645-680}
\medskip
\item{[W]} {Weyl, H., {\it  \"Uber die Gleichverteilung von Zahlen mod Eins\/}, Math. Ann., {\bf 77}, 1956, pp. 313--352 } 
\medskip

\enddocument

\ref 
\key B
\by  Billingsley, P.
\book Convergence of Probability Measures, 2nd Edition
\publ John Wiley
\yr 1999
\endref

\ref \key Bl
\by  Blaschke, W.
\book Vorlesungen \"uber Integralgeometrie,  Bd.I 1936; Bd.II 
\publ Teubner Verlag, Leipzig
\yr 1937
\endref

\ref \key B-C
\by  Bailey, C.H. and Crandal, R.
\paper  Random Generators and Normal Numbers
\jour Experiment. Math.
\publ A.K. Peters Ltd.
\vol 11 (4) \yr 2002 \pages 527--546
\endref

\ref \key B-T
\by  Bott, R. and Tu, L.
\book Differential Forms in Algebraic Topology
\publ Springer Verlag
\yr 1982
\endref

\ref \key C
\by Chern, S.S.
\paper  On Integral Geometry in Klein Spaces 
\jour Ann. of Math.
\vol 43 \yr 1942  \pages 178--189 
\endref

\ref \key Cr
\by Crofton, M.W.
\paper  Probability
\jour Encyclopedia Britannica, 9th Ed.
\vol 19 \yr 1885 \pages 768--788
\endref

\ref \key E
\by L'Ecuyer, P.
\paper Uniform Random Number Generators: A Review
\jour Proceedings of the I997 Winter Simulation Conference
\vol  \yr 1997 \pages 127--134
\endref

\ref \key F
\by Franklin, J.N.
\paper Deterministic Simulation of Random Processes
\jour Math. Computation
\vol 17 \yr 1963\pages 28--57 
\endref

\ref \key Fe
\by Federer, H.
\paper Colloquium Lectures on Geometric Measure Theory
\jour Bull. Amer. Math. Soc.
\vol 84 \yr 1978\pages 291--338 
\endref

\ref \key F-W1
\by Fuselier, E. and Wright,G.
\paper A High-Order Kernel Method for Diffusion and Reaction-Diffusion Equations on Surfaces
\jour Journal of Scientific Computing
\vol 56 (3) \yr 2013\pages 535--565 
\endref

\ref \key F-W2
\by Fuselier, E. and Wright,G.
\paper Scattered Data Interpolation on Embedded Submanifolds with Restricted 
Positive Definite Kernels: Sobolev Error Estimates
\jour SIAM Journal on Numerical Analysis
\vol 50 (3) \yr 2012\pages 1753--1776 
\endref

\ref \key H
\by  Hirsch, M.W.
\book Differential Topology
\publ Springer Verlag
\yr 1976
\endref

\ref \key J
\by  Judd, K.L.
\book Numerical Methods in Economics
\publ MIT Press
\yr 1998
\endref

\ref \key Ka
\by Kantorovich, L. V.
\paper  On NewtonÕs method for functional equations 
\jour Dokl. Akad. Nauk. SSSR, 
\vol 59 \yr 1949  \pages 1237--1240 
\endref

\ref \key K1
\by  Knuth, D.
\book The Art of Computer Programming, Volume 2:   Seminumerical Algorithms
\publ Addison-Wesley
\yr 1997
\endref

\ref \key K2
\by  Knuth, D.
\paper Construction of a Random Sequence
\jour BIT
\vol 5 \yr 1965 \pages 246--250 
\endref

\ref \key {K-R}
\by Klain,D. and Rota G.-C.
\paper  Introduction to Geometric Probability
\jour  Cambridge Univ. Press, Cambridge,
\vol \yr 1997 
\endref

\ref \key K-N
\by  Kuipers, L. and Niederreiter, H
\book Uniform Distribution of Sequences
\publ   Wiley \publaddr New York
\yr 1974
\endref

\ref \key L-Z
\by Liang, J. and Zhao, H.
\paper  Solving partial differential equations on point clouds, 
\jour SIAM Journal on Scientific Computing
\vol 35(3) \yr 2013  \pages 1461--1486 
\endref

\ref \key {M-N}
\by Matsumoto, M.; Nishimura, T. 
\paper  Mersenne twister: a 623-dimensionally equidistributed 
     uniform pseudo-random number generator
\jour ACM Transactions on Modeling and Computer Simulation 
\vol 8 \yr 1998  \pages 3--30 
\endref

\ref \key N
\by Nash, J.
\paper The imbedding problem for Riemannian manifolds
\jour Ann. of Math. 
\vol 63 (1) \yr 1956 \pages 20--63   
\endref

\ref \key O
\by Ortega, J.M. 
\paper The Newton-Kantorovich Theorem
\jour The American Math. Monthly
\vol 75 (6)  \yr 1968\pages 658--660
\endref

\ref \key P
\by Palais, R.
\paper The Visualization of Mathematics: Towards a Mathematical Exploratorium
\jour Notices of the AMS
\vol 46 (6) \yr 1999\pages 647--658
\endref

\ref \key PF
\by Paiva, J.C. \'Alvarez and  Fernandes, E.
\paper Gelfand Transforms and Crofton Formulas
\jour Sel. Math., New Ser.
\vol 13  \yr 2007\pages 369--390
\endref

\ref \key {R}
\by Riesz, F.
\paper  Sur les op\'erations fonctionnelles lin\'eaires,
\jour C.R. Acad. Sci., Paris 
\vol 149 \yr 1909  \pages 974--977 
\endref

\ref \key {S1}
\by Santalo, L.
\paper  Introduction to integral geometry
\jour  Hermann, Paris
\vol \yr 1953  \pages 36-88 
\endref

\ref \key {S2}
\by Santalo, L.
\paper  Integral Geometry and Geometric Probability
\jour  Addison-Wesley
\vol \yr 1976  
\endref

\ref \key {Sp1}
\by Spivak, M.
\paper  Calculus on Manifolds
\jour  W.A. Benjamin
\vol \yr 1965 
\endref

\ref \key {Sp2}
\by Spivak, M.
\paper  A Comprehensive Introduction to Differential Geometry Vol. I
\jour  Publish or Perish
\vol \yr 1999
\endref

\ref \key {U}
\by Metropolis, N. and Ulam, S. 
\paper  The Monte Carlo method
\jour Journal of the American Statistical Association
\vol 44 \yr 1949  \pages 335--341 
\endref

\ref \key {VN}
\by von Neumann, J.
\paper  Various techniques used in connection with random digits, Monte Carlo Method
\jour  Nat. Bur. Stand., Appl. Math. Series
\vol 12 \yr 1951  \pages 36-88 
\endref

\ref \key W
\by Weyl, H.
\paper \"Uber die Gleichverteilung von Zahlen mod. Eins
\jour Math. Ann.
\vol 77 \yr 1916 \pages 313--352 
\endref

Marsaglia, G. (1972). "Choosing a Point from the Surface of a Sphere". 
Ann. Math. Stat. 43 (2): 645Ð646. doi:10.1214/aoms/1177692644.


\vfil\eject   

\remark{Conjecture} 
If $\{S_j\}$ is an equidistributed sequence of $k$-planes in $\Reals^n$ and $\{T_j\}$ is an equidistributed sequence 
of $l$-planes in $\Reals^n$  with $S_j \perp T_j$, the $\{S_j \oplus T_j\}$ is an equidistributed sequence of 
$k+l$ planes in $\Reals^n$.

\endremark

\centerline {Density Variation for the Early Point Cloud Algorithm}

\medskip

We will show that with the early point cloud algorithm, in which a surface is intersected with 
three orthogonal families of lines, the local density of points of the resulting point cloud changes 
by a factor of $\sqrt 3$ as the normal to the surface varies.

  Let $e \in \S^2$ and let $\Cal F_e$ be a family of lines with direction $e$ that intersect
the plane $e^\perp$ with a density of $N$ points per unit area. (Equivalently, let $\Cal P$ 
be a uniformly distributed family of points of $e^\perp$ with a density of $N$ points per unit area,
and let $\Cal F_e$ be the lines $ t \mapsto t e + p$.)

  If $\nu$ is another unit vector in $\Reals^3$ and $\Sigma$ is a small open set in $\nu^\perp$, 
then the area of the projection of $\Sigma$ onto $e^\perp$ is $| e\cdot\nu|$ times the area of 
$\Sigma$, and it follows that the density of the set of points where the lines of $\Cal F_e$ meets 
$\Sigma$ is $N | e\cdot\nu|$.

  So, if $e_1, e_2, e_3$ is an orthonormal basis for $\Reals^3$ and $\Cal F_{e_i}$ is a family of 
lines with direction $e_i$ that intersect the plane $e_i^\perp$ with a density of $N$ points per 
unit area, then the density of the set of points where the lines of $\Cal F := \Cal F_1 \cup \Cal F_2 \cup \Cal F_3$ 
meets $\Sigma$ is $\rho(\nu,N) := N (| e_1 \cdot\nu| +  | e_2 \cdot\nu| + | e_3 \cdot\nu|)$.

  The minimum value, $N$, of $\rho(\nu,N)$ is assumed when $\nu =\pm e_i$ and its maximum 
value, $N \sqrt 3$, is assumed when $\nu = {\sqrt 3\over 3} (\pm e_1  \pm e_2 \pm e_3)$, and we 
note that the ratio $\sqrt 3$ of the maximum to the  minimum value is independent of $N$.

   Of course at a point $p$ of a smooth surface $\Sigma \subset \Reals^3$ where the 
normal to $\Sigma$ is $\nu$, the local density at $p$ of the point cloud formed by intersecting 
$\Sigma$ with $\Cal F$ will be $\rho(\nu,N)$, so this local density will also vary between $N$ 
and $N \sqrt 3$, and the max to min ratio is again independent of $N$ and for a convex surface 
will be $\sqrt 3$.

  An interesting experiment would be to construct the point cloud for the pyramid with the 
  four vertices $O, e_1, e_2, e_3$. The three faces parallel to a coordinate plane would all 
  have density $N$ while the fourth face, spanned by the $e_i$, would have density $N \sqrt 3$ ,
  and rotating this pyramid should show the difference pretty clearly. Another interesting 
  experiment would be to rotate $\S^2$ along the great circle joining say $e_1$ and 
  ${\sqrt 3\over 3} (\pm e_1  \pm e_2 \pm e_3)$, and watching how the density changes.

\enddocument

  If $f: \Sigma^k \to \Reals$ is a bounded continuous function, then by the definition of 
$f^\Sigma$ in section 8, it is immediate from the above algorithm that:

{\narrower\narrower
\remark{Observation} If $N := N_1 N_2 \cdots N_m$ then the first $m$ lines of the sequence $\{\ell_j\}$ 
meet $\Sigma$ in the $N$ points $x_1, \ldots, x_N$. Hence if $f: \Sigma^k \to \Reals$ is a bounded 
continuous function, then by the definition of $f^\Sigma$ in section 8, 
$$   {1 \over N} \sum_{1\le i_1,\ldots,i_k \le m}f^\Sigma(\ell_{i_1},\ell_{i_2},\ldots,\ell_{i_k})$$
\endremark
 }

\smallskip\noindent
We will use this observation to prove that the sequence

\ref \key P
\by Palais, R.
\paper The Visualization of Mathematics: Towards a Mathematical Exploratorium
\jour Notices of the AMS
\vol 46 (6) \yr 1999\pages 647--658
\endref

The Visualization of Mathematics: Towards a Mathematical Exploratorium
Palais, R. 
NOTICES OF THE AMS
VOLUME 46, NUMBER 6
1999

\medskip\noindent
 We associate to $\Rnd_\Sigma$ an increasing sequence of ``point clouds'', 
$PC_N(\Sigma)$, namely the $N$-element subsets of $\Sigma$ consisting of the first 
 $N$ elements of the sequence $\{x_n\}$. As $N$ becomes large, these point clouds 
 can be used for many purposes as a good approximation to the surface $\Sigma$, 
 and we will discuss below various senses in which this is true. 

  As our notation suggests, we will think of the sequence of points $\{ \Rnd_\Sigma(n) \}$  
as being in some sense ``randomly distributed in $\Sigma$'', much as a so-called   
random number generators (RNG) is usually thought of as a sequence of 
points randomly distributed in $[0,1)$. But, as is well-known, the problem of giving 
a satisfactory definition of randomness for a deterministic algorithm is more difficult 
than one might at first expect, and we will put off to the Appendix the discussion of how 
to make this precise. 

The algorithms $\Rnd_\Sigma$ that we describe will be built from a given random 
number generator $\Rnd_{[0,1)}$,  for which  we assume at a minimum that  the
sequence $\Rnd_{[0,1)}(n)$ is uniformly distributed in $[0,1)$. This property of a RNG 
is also often referred to as being $1$-distributed, and there is also a sequence of 
increasingly stronger properties for an RNG, called being {\it $k$-distributed\/}. 
Namely, $\Rnd_{[0,1)}$ is $k$-distributed if the sequence of points in 
$[0,1)^k$  given by:
$$\Rnd_{[0,1)^k}(n) := (\Rnd_{[0,1)}(n), \Rnd_{[0,1)}(n+1), \ldots, \Rnd_{[0,1)}(n+k))$$
is uniformly distributed in $[0,1]^k$,  and we say that a RNG is {\it $\infty$-distributed\/} 
if it is \hbox{$k$-distributed} for all $k$. (See \cite{K}  section 3.5 for a full discussion of these 
properties.)   Randomness for a RNG  is sometimes defined 
as being \hbox{$\infty$-distributed}, and it is an analogous condition that we will use as 
our definition of a point cloud being pseudo-random in the appendix.  We will see 
that if $\Rnd_{[0,1)}$ is pseudo-random in this sense then so are the algorithms 
$\Rnd_\Sigma$ that we will build using it.

  Before describing our methods for constructing $\Rnd_\Sigma$ for various types of 
surfaces $\Sigma$, we first recall the definitions and some properties of these types.

\section{Appendix: When is a Point Cloud Random?}

  \subsection{Notation}
  
    In what follows $X$ will denote a compact metric space, and we recall 
 that a {\it Radon measure\/}  on $X$ is a finite measure on the Borel sets of $X$, 
(the $\sigma$-algebra generated by the open sets)  that is inner-regular, 
i.e.,  the measure of any Borel set $B$ is the sup of the measures of compact 
subsets of $B$.  We recall next an important  characterization of Radon measures
that was made into a alternative definition by Bourbaki. \par

    We denote by $C(X)$ the space of continuous, real-valued functions on $X$
 with the  sup norm, $||f||_\infty = \sup|\{ f(x) \mid x \in X\}$. Clearly, integration with
 respect to a Radon measure $\mu$ defines a continuous linear functional $\ell_\mu$ 
 on $C(X)$  that is  positive, i.e., non-negative on positive functions.  The Riesz 
 Representation Theorem [R] states that, conversely,  every positive continuous linear 
 functional $\ell$ on $C(X)$ is of the form $\ell_\mu$ for a unique 
 Radon measure. 
 
  We will denote by $\mu_X$  some fixed non-trivial Radon measure on $X$ 
and  for each positive integer $k$ we let $\mu_{X^k}$ denote the associated
Fubini  product measure on $X^k$;  $\mu^*_X$ is the associated probability 
measure on $X$,  obtained by dividing by the measure of the whole space, 
and similarly for  $\mu^*_{X^k}$.
 
  Let $\{x_n\}$ be a sequence of points of $X$ and for each $n$ let $S(n)$  be a 
statement about $x_n$. If $\nu(N)$ is the number of statements 
$S(1), S(2),\ldots S(N)$ that are true, we will write $Pr(S(n)) = \lambda$ 
and say that $S(n)$ holds with probability $\lambda$ if 
$\lim _{N\to \infty} {\nu(N) \over N} = \lambda$. \par
   
    $\Rnd_X$ will be an algorithm that associates to each positive integer $n$ 
 a point $x_n = \Rnd_X(n) \in X$, and $PC_N(X)$ will denote the set 
$\{x_1,x_2, \ldots , x_N\}$ consisting of the first $N$ elements of the 
sequence $\{ x_n \}$. We call $\{PC_N(X)\}$ the sequence of {\it point clouds\/} 
generated by $\Rnd_X$.  \par

 For each positive integer $k$ there is an algorithm $\Rnd_{X^k}$ naturally 
 associated to $\Rnd_X$ for producing a sequence of points in $X^k$; namely
 $$\Rnd_{X^k}(n) := (x_n, x_{n+1},\ldots,x_{n+k-1}).$$
 
 \goodbreak

\subsection{Definitions of $k$-Distributed and  $\infty$-Distributed or Pseudo-Random}

\definition {Definition}

The sequence  $\{x_n\}$ is said to be  {\it uniformly distributed in $X$ \/} with respect to 
 $\mu_X$ if  for every open subset $O$ of $X$ the probability 
 $Pr(x_n \in O)$ that $x_n$  is in $O$ is equal to $\mu^*_X(O)$, and $\{x_n\}$ is said to be 
{\it $k$-distributed in $X$ \/} with respect to the measure $\mu_X$ if the sequence
$\Rnd_{X^k}(n)$ defined above  is uniformly distributed in $X^k$ with respect to 
the measure $\mu_{X^k}$.  Finally,  $\{x_n\}$  is said to be  {\it $\infty$-distributed in $X$ \/}  
or {\it pseudo-random\/} if it is $k$-distributed in $X$ for every positive integer $k$. 
\enddefinition

\definition {Definition} 
An integrable real-valued function $f:X^k \to \Reals$ is said to be {\it Monte Carlo Integrable\/} 
with respect to the sequence  $x_n = \Rnd_X(n)$ if 
$$  
   {1\over \mu_X(X)^k}    \int_{X^k} f(x_1,\ldots,x_k) \, d\mu_{X^k} = 
           \lim_{n \to \infty} {1\over n} \sum_{j =1}^n f(x_j,\ldots,x_{j+k-1} )
$$

\enddefinition

  The following remarks are obvious:
  
\item{i)} The sequence $x_n$ is $k$-distributed if and only if the characteristic function 
of every open subset of $X^k$ is MC-Integrable. 

\item{ii)} The set of MC-Integrable.functions is a linear subspace of $L^1(X^k,\mu_{X^k})$.

\section{Code}

Generating point clouds on triangulated or parametric surfaces.

I. The surface is decomposed as a set of triangles (sets of three vertices in three-dimensional
space) or as a set of patches, typically rectangles, in the domain in two-dimensional space
that is mapped to three-dimensional space by functions parametrizing the surface, x(u,v),y(u,v),z(u,v).

Examples of triangulated surfaces: 

A polyhedron whose faces all are triangles, such as a tetrahedron, an octahedron. 
A polyhedron whose non-triangular faces are decomposed into triangles by the addition of edges. 
A sphere or other closed surface that is approximated by a polyhedron, then triangulated
as above. There are many known algorithms for triangulating a surface.

Examples of parametrized surfaces decomposed into rectangular patches: 

The standard torus, the unit sphere, the ellipsoid, a section of the paraboloid, a section of the hyperboloid.

The standard torus: The (u,v) rectangular domain in the plane $ [0,2\pi] \times [0,2\pi]$ is mapped to 
$(x,y,z)=(\cos u \cos v, \cos u \sin v, \sin u)$ \par

\medskip\noindent
II. Distribute the points among triangles in a triangulation or  (typically rectangular) grid patches on a surface parametrized by functions, pseudo-randomly and in a manner that is asymptotically proportional to the area represented by each triangle or patch.

 A.  Initialize the cumulative area distrubution function.
 
  Approximate areas of each triangle or patch. 
  
  If the surface is is a discrete triangulation, the function Area(patch) approximates the area of a patch using the formula for the area of a three-dimensional triangle in terms of its vertices, one-half of the magnitude of the cross-product of two edges (where an edge is the difference of two vertices).

If the surface is given analytically, Area(patch) approximates the area of a patch using the magnitude of the determinant of the cross-product of the partial derivatives of (x,y,z) with respect to u and v at a vertex.
  
  Store the partial sums of these area in an array we call range[.]. So, range[0] is the approximate area of the first triangle or patch, range[1], is the sum of the approximate areas of the first two triangles or patches, and if there are N triangles or patches, range[N-1] is the approximate area of the surface. If we consider a `virtual' array element, ranges[-1]=0, then neighboring values in the array give the endpoints of a sequence of consecutive non-overlapping intervals [ range[i-1] , range[i] ) whose length is the approximate area of the ith triangle or patch, and (range[i]-range[i-i])/range[N-1]  gives the proportion of the total area contributed by each triangle or patch.
  
  The pseudocode for initializing the cumulative area distribution function is:

\verbatim
Area_approx=0; //  initialize cumulative area function
for (patch=0;patch<patches;patch++){
Area_approx+=Area(patch);  // Add patch area to cumulative area
range[patch]=Area_approx; //Store result in range[] array
}
|endverbatim

After implementing this step,  
\verbatim Area_approx=range[N-1]  |endverbatim
 is the approximate total area of the surface.

 B.  Assign a point to a triangle or patch based on the cumulative area function

 Let rand be a number taken from uniformly distributed pseudo-random numbers on [0,1] then multiplied by the total area, 
\verbatim Area_approx, |endverbatim
 so rand now lies in 
\verbatim  [0,Area_approx].  |endverbatim
 Then the likelihood that such a number is in the interval corresponding to the ith patch is equal to the proportion of area it contributes to the total area of the surface.

Use bisection to find the interval in which rand lies in an efficient manner.

Initially, and inductively, at each step, we keep rand in the interval [ range[left], range[right] ), inclusive on the left, exclusive on the right. The left and right indices are always integers so we can bisect or approximately bisect the remaining sub-array into two consecutive non-overlapping sub-arrays whose lengths differs by no more than 1. The width, right-left is strictly decreasing by bisection to 1. 

For example if the difference right-left is even, e.g., right=8 and left=2, whose difference is 6, and whose midpoint is 5, we bisect  exactly and determine if rand is in  [ range[2],range[5] ) or [ range[5],range[8] ).
If the difference right-left is odd, e.g., right=7 and left=2, whose difference is 5, and whose midpoint is 4.5, we bisect  approximately into consecutive non-overlapping sub-arrays whose length differs by 1
and determine if rand is in  [ range[2],range[4] ) or [ range[4],range[7] ) . 

In this implementation, the subarray  that is smaller is on the same side as the interval in which rand was found at the previous stage, for example, if it was found in the interval corresponding to the left of the previous two subintervals, that left subarray is split into two with the left being smaller or equal to the right, and if it was in the interval corresponding to the right, that right subarray is split into two with the right being smaller or equal to the left.

 The pseudocode for scaling rand and finding the interval in which the scaled value lies is:

\verbatim
function findinterval(rand){
rand=Area_approx*Math.random(); 
while (halfwidth>=1) {  
	if (rand>=range[mid]) 
	{ left=mid; halfwidth=Math.floor((right-left)/2); 
	   mid+=halfwidth; } 
	            // preserves rand in [range[left], range[right]) . 
	            //Since halfwidth<= half the width, (left,mid) 
	            //ends up smaller or equal to (mid,right)
	else 
	{ right = mid; halfwidth=Math.floor((right-left)/2);
	    mid-=halfwidth; }
	            // preserves rand in [range[left], range[right]) . 
	            //Since halfwidth<= half the width, (mid,right) 
	            //ends up smaller or equal to (left,mid)
		}
	return right; // rand is in [range[right-1], range[right])
}
|endverbatim

III. Distribute points within a triangle in a triangulation, or a (typically rectangular) grid patch of a surface parametrized by functions.

Once a point is assigned to a patch or triangle, assign it to a random point in that triangle or patch on the surface.

For a (typically rectangular) grid patch of a surface parametrized by functions:

 Let randu and randv be two numbers taken from uniformly distributed pseudo-random numbers on [0,1]. Scale them to the grid scales, [0,0]x[du,dv]. Translate to the particular patch in the domain. Map them to the corresponding patch on the surface using the parametrization functions.

The pseudocode for doing this is:

\verbatim
patch=findinterval(rand);
	i=Math.floor(patch/num_v);
	j=patch-i*num_v;
	randu=uleft+(i+Math.random())*du;
	randv=vlower+(j+Math.random())*dv;
	xrand=x(randu,randv);
	yrand=y(randu,randv)
	zrand=z(randu,randv);
	
|endverbatim

For a triangle, generate uniformly distributed barycentric coordinates on the standard equilateral
simplex, x+y+z=1, x,y,z>=0 by acceptance/rejection, and map them to the assigned triangle
by forming a convex combination of its vertices using these barycentric coordinates as coefficients.

The pseudocode for doing this is:
\smallskip

\verbatim
triangle=findinterval(rand);
	r1=Math.random(); 
	r2=Math.random();
	if ((r1+r2)<1){
	r3=1-(r1+r2);
	randompoint[0]=r1;
	randompoint[1]=r2;
	randompoint[2]=r3;
				 }
	 else { 
	randomtriangle( randompoint );
	 } // try again
	
    // The jth compontent is a convex combination of
    // jth component of the three vertices i=0,1,2 of
    // the chosen triangle
	for (i=0;i<3;i++){
	for (j=0;j<3;j++){
randomdot[j]+=randompoint[i]* vertex [triangle] [i] [j] ;
 }
   }
|endverbatim

\verbatim
int RotationMatrix(double uI[3], double uF[3], double T[3][3])
{

/*
 * A routine that takes as input an `initial' unit vector, 
 * uI[3], and a `final' unit vector, uF[3], not equal to -uI, 
 * and returns the unique rotation matrix R[3][3] satisfying
 *  R uI = uF and R v = v if v is orthogonal to uI and uF.
 * Author: Bob Palais
 * "A New Formula For Computing A Rotation Matrix"
 */
 
 #define TOL 1E-16
 int nonzero_sum = false;
 double s_dot_s, s0,s1,s2;
 double out00,out01,out02,out10,out11,out12,out20,out21,out22;
 double uFF0, uFF1, uFF2, uI0, uI1, uI2;
 double k, ks0, ks1, ks2, ks00, ks11, ks22, ks01, ks12, ks20;

 uFF0=uF[0]; uFF1=uF[1]; uFF2=uF[2];
 uI0=uI[0]; uI1=uI[1]; uI2=uI[2];

 // compute s=uI+uF and 2uF
 s0=uFF0+uI0; s1=uFF1+uI1; s2=uFF2+uI2;

 // compute  2uF
 uFF0=uFF0+uFF0; uFF1=uFF1+uFF1; uFF2=uFF2+uFF2;

 // compute diagonal of outer product out = 2uF uI^T
 out00=uFF0*uI0; out11=uFF1*uI1; out22=uFF2*uI2;
 // 3 multiplications

 // use it to compute  (s,s) = 2 + (uI,2uF) 
 s_dot_s = 2.0 + out00 + out11 + out22;

 if (s_dot_s > TOL ) // check for uF = -uI
 { 
  nonzero_sum = true;
  k = 2.0 / s_dot_s; // compute normalization factor k=2/(s,s)
  // 1 division

  // 2uF uI^T, off-diagonal components
  out01 = uFF0*uI[1]; out12 = uFF1*uI[2]; out20 = uFF2*uI[0];
  out10 = uFF1*uI[0]; out21 = uFF2*uI[1]; out02 = uFF0*uI[2];
  // 6 multiplications

  //  compute components of (ks) and (ks)s^T
  ks0 = k*s0; ks00 = ks0*s0; ks01 = ks0*s1;  
  ks1 = k*s1; ks11 = ks1*s1; ks12 = ks1*s2;
  ks2 = k*s2; ks22 = ks2*s2; ks20 = ks2*s0;
  // 3 multiplication +  6 multiplications

  // now build the matrix
  T[0][0] = 1.0 + out00  - ks00;
  T[1][1] = 1.0 + out11  - ks11;
  T[2][2] = 1.0 + out22  - ks22;
  T[0][1] = out01 - ks01;
  T[1][2] = out12 - ks12;
  T[2][0] = out20 - ks20;
  T[1][0] = out10 - ks01;
  T[2][1] = out21 - ks12;
  T[0][2] = out02 - ks20;
 }
 return nonzero_sum;
}
|endverbatim

\subsection{Uniformly distributed pseudo-random dots on a parametric surface}

Using the parametric definition of the
surface, $X(u,v), Y(u,v), Z(u,v)$, we can
distribute dots pseudo-randomly on a parametrically defined surface.
First we first compute the area $A_j, ~ \j=1,2,\dots$ of the
image of each patch $P_j$ from the magnitude
of the cross product of two adjacent edges. We use these to
construct a sequence of points $x_j$
such that $x_0=0$ and $x_j=x_{j-1}+A_j, ~ j=1,2,\ldots$.
This creates a subdivision of the interval $[0,A)$ where
$A=\sum_j A_j = \sum (x_j-x_{j-1} ) $.
Then we can use pseudo-random numbers $r_i$ from a distribution 
scaled to be $\infty$-distributed on the interval $[0,A)$,  
to distribute dots to patches in proportion to their area by
assigning the $i$th dot to the $j$th patch if $r_i \in [ x_{j-1}, x_j )$.
Two more pseudo-random numbers are taken from
a distribution scaled to be uniform in
the standard rectangular patch $[0,du] \times [0, dv]$ in the domain.
This resulting point is translated to the patch $P_j$ and mapped to the
range.

(We can accelerate the process of finding the interval in
which $r_i$ is found using...a binary tree?)

\subsection{Implementing the Cauchy-Crofton measure}

To implement this algorithm, we need to generate uniformly
distributed vectors on the unit sphere $\u_j \in S^2$ in $\Reals^3$,
and uniformly distributed points $p_j$ on the disc $D_j$ of radius $R$
orthogonal to $\u_j$, centered at the origin in $\Reals^3$.
Then the lines $l_j (s)= p_j+ s \u_j=(x(s),y(s),z(s))$ will be uniformly 
distributed and our uniformly distributed dots on the implicitly defined 
surface $S=\{(x,y,z) | f(x,y,z)=0\}$
are obtained by finding intersections of $l_j(s)$ with $S$ using
familiar methods from ray tracing to determine where $f(x(s),y(s),z(s))$
changes sign.

Perhaps the simplest and most efficient algorithm to 
 obtain $\u_j$ distributed uniformly with respect to area on
$S^2 =\{ \u \in \Reals^3 | ||\u||=1 \}$ is to generate
uniformly distributed points $\r_j$ in the cube $[-1,1]^3$ and
exclude those such that $|| \r_j^2||>1$, which involves computing
a sum of squares. Approximately  ${4 \over 3}\pi/8 \approx 0.52$ of 
the candidate points are accepted, and the resulting points are
 uniformly and isotropically distributed in the unit ball in $ \Reals^3$.
We normalize them using one square root, one division, and 
three multiplications per point, to obtain unit vectors distributed
uniformly with respect to area on $S^2$.

We can then adapt the two-dimensional version of the first steps
of this process to uniformly distributed points on $D_j$.
We first generate
uniformly distributed points $\r_j$ in the square $[-1,1]^2$ in the
$x-y-$plane and exclude those such that $|| \r_j^2||>1$.
Approximately  $\pi/4 \approx 0.78$ of the candidate points are
accepted, and the resulting points are uniformly distributed in
the unit disc in the $x-y-$plane. Next we scale $r_j$ by a factor of
$R$ to the disc of radius $R$, and call the result $\v_j=R\r_j$. 
The last step is to rotate $\v_j$ to the disc $D_j$. We accomplish
this efficiently by performing two reflections that implement the
rotation $\R$ taking the vector $\e_3=<0,0,1>$
orthogonal to our standard disc to the vector $\u_j$ orthogonal
to $D_j$, and fixes the axis orthogonal to both. Specifically,
we let $\s=\u_j+\e_j$ and $\w=2(\v_j \cdot \e_3)\e_3-\v_j)$, 
then $\R\v_j=2(\w_j \cdot \s)\s-\w_j$. See [TransvectionForRotations]
for more details and other applications of this technique.

Several other methods can be used to implement each of the
steps above. Instead of restricting uniformly distributed
points in the cube to the ball and normalizing, we can
obtain isotropically distributed points in $\Reals^3$
by taking each of their components from a normally distributed
random variable, then normalizing and proceeding as before.
This approach works in arbitrary dimension and goes
back to Maxwell and Boltzmann [Reference?]
is well-known in statistical mechanics for
obtaining Gibbs measure on the spherical shell of
states in phase space having energies in a given interval.
\footnote{
A related though more involved approach has also been
suggested, in which a point follows a random walk starting
from the origin and $\u_j$ is the point at which it first escapes
the unit ball.}

We may also obtain uniformly distributed points on the unit
sphere using spherical coordinates, though the fact that
the co-latitude $\phi$ and longitude $\theta$  cannot be individually uniformly
distributed is the subject of the well-known Borel's paradox, 
discussed in [Billingsley]. To address this, we wish points
to be distributed proportionally to the area element $\sin \phi d\theta d\phi$.
(The circumference of a latitude circle at $\phi=\phi_0$ is
proportional to $\sin \phi_0$.) We may incorporate this distribution
using an acceptance-rejection method like the one used above.
First we choose a candidate co-latitude
$\phi_j$ from a uniform distribution on $[0,\pi)$, then
choose a number $r_j$ from a uniform distribution on $[0,1)$.
If $r_j< \sin \phi_j <1$ we  scale it to $\theta_j= 2\pi r_j/\sin \phi_j $ in the
interval $[0,2\pi)$, and accept the pair $(\theta_j, \phi_j)$ as the spherical coordinates
of the desired $\u_j = < \sin \phi \cos \theta , \sin \phi \sin \theta , \cos \phi >$.
If $r_j > \sin \phi_j <1$ we reject it and search for a point with a different $\phi_j$.
If $\phi=\pi/2$, any value of $\r_j$ will give a point on the equator, and if $\phi = 0$
or $\pi$, no values of $r_j$ will be accepted to give points at the poles.

The spherical coordinates approach can also be implemented by inverse
transform sampling.  If we rewrite the spherical area measure on
the sphere $\sin \phi d\theta d\phi = d\theta d(-\cos \phi )$
then we can invert using $\arccos(1-2x)$ to map uniformly distributed 
$x_j \in [0,1)$ to $\phi_j \in [0,\pi)$ with the distribution $d(-\cos \phi )=\sin \phi$.

    * Patrick Billingsley (1979). Probability and Measure. New York, Toronto, London: 
    John Wiley and Sons. 

\enddocument